\documentclass[12pt]{article}
\usepackage{mathrsfs}
\usepackage{latexsym, amssymb, amsmath, amscd, amsfonts, epsfig, graphicx, colordvi,verbatim,ifpdf}
\usepackage{amsfonts, amsmath, amssymb}
\usepackage{amssymb,amsfonts,amsmath,latexsym,epsfig,cite, psfrag,eepic,color}
\usepackage{amscd,graphics}
\usepackage{latexsym, amssymb,  amsmath,amscd, amsfonts, epsfig, graphicx, colordvi,amsthm}

\usepackage{graphicx}
\usepackage{color}
\usepackage{ifpdf}
\usepackage{fancybox}

\usepackage{float}

\newtheorem{thm}{Theorem}[section]

\newtheorem{lem}[thm]{Lemma}

\def\pf{\noindent{\it Proof.} }
\def\qed{\nopagebreak\hfill{\rule{4pt}{7pt}}
\medbreak}

\numberwithin{equation}{section}

\def\qed{\nopagebreak\hfill{\rule{4pt}{7pt}}
\medbreak}

\newlength{\boxedparwidth}
  {\begin{center} \begin{tabular}{|@{\hspace{.315in}}c@{\hspace{.15in}}|}
                  \hline \\ \begin{minipage}[t]{\boxedparwidth}
                  \setlength{\parindent}{.25in}}%
  {\end{minipage} \\ \\ \hline \end{tabular} \end{center}}
\allowdisplaybreaks
\begin{document}
\begin{center}

 {\Large \bf $M_2$-Ranks of overpartitions modulo $6$ and $10$}
\end{center}
\vskip 0.2cm

\begin{center}
{Helen W.J. Zhang} \vskip 3mm

  Center for Applied Mathematics\\[2pt]
Tianjin University, Tianjin 300072, P.R. China\\[5pt]

  wenjingzhang@tju.edu.cn

\end{center}

\vskip 6mm \noindent {\bf Abstract.} In this paper, we obtain inequalities on
$M_2$-ranks of overpartitions modulo $6$. Let $\overline{N}_2(s,m,n)$ to be the number of overpartitions of $n$ whose $M_2$-rank is congruent to $s$ modulo $m$. For $M_2$-ranks modulo $3$, Lovejoy and Osburn derived the generating function of  $\overline{N}_2(s,3,n)-\overline{N}_2(t,3,n)$,
which implies the inequalities $\overline{N}_2(0,3,n)\geq\overline{N}_2(1,3,n)$. For $\ell=6, 10$, we consider
the generating function
$\overline{R}_{s,t}(d,\ell)$ of the $M_2$-rank differences $\overline{N}_2(s,\ell,\ell n/2+d)+\overline{N}_2(s+1,\ell,\ell n/2+d)-\overline{N}_2(t,\ell,\ell n/2+d)-\overline{N}_2(t+1,\ell,\ell n/2+d)$.
By the method of Lovejoy and Osburn, we derive a formula for $\overline{R}_{0,2}(d,6)$.
This leads to the inequalities for $n\geq0$,
$\overline{N}_2(0,6,3n)\geq\overline{N}_2(2,6,3n)$ and $\overline{N}_2(0,6,3n+1)\geq\overline{N}_2(2,6,3n+1)$.
Based on the valence formula for modular functions, we compute $\overline{R}_{0,4}(d,10)$ and $\overline{R}_{1,3}(d,10)$.
In particular, we notice that the generating function $\overline{R}_{0,2}(2,6)$ can be expressed in terms of the third order mock theta function $\rho(q)$, and the generating functions $\overline{R}_{0,4}(4,10)$, $\overline{R}_{1,3}(1,10)$ and $\overline{R}_{1,3}(4,10)$ can also be expressed in terms of the tenth order mock theta functions $\phi(q)$ and $\psi(q)$. \vskip 0.3cm

 \vskip 0.3cm

 \vskip 0.3cm

\section{Introduction}
The rank of a partition was introduced by Dyson \cite{Dyson-1944}, which equals to the largest part minus the number of parts.
%Responding to the conjecture of Dyson,
%Dyson's rank of partitions gives combinatorial interpretations for Ramanujan¡¯s congruences for partition
%function modulo 5 and 7
Let $N(s,\ell,n)$ denote the number of partitions of $n$ with rank congruent to $s$ modulo $\ell$.
Atkin and  Swinnerton-Dyer \cite{Atkin-Swinnerton-Dyer-1954} obtained  generating functions for rank differences $N(s, \ell, \ell n + d)- N(t, \ell, \ell n + d)$ with $\ell = 5$ or $7$ and $0\leq d, s, t < \ell$, which lead to combinatorial interpretations of Ramanujan's congruences modulo $5$ and $7$.
The generating functions for the rank difference $N(s, \ell, \ell n + d)- N(t, \ell, \ell n + d)$  with $\ell=11$ have  also  been determined by Atkin and Hussain \cite{Atkin-Hussain-1958}, which implies Ramanujan's congruence modulo $11$.
Since then, the ranks of partitions with moduli other than $5$, $7$ and $11$ have been extensively studied, see, for example, \cite{Lewis-1997,Lewis-1991,Gadea-1992,Mao-2013}.

%The rank of a partition is now understood to be a special case of a more general rank which is defined on an overpartition.
As to overpartitions, there are two kinds of ranks defined by Lovejoy:  $D$-rank \cite{Lovejoy-2005} and $M_2$-rank \cite{Lovejoy-2008}. An overpartition \cite{Corteel-Lovejoy-2004} of a nonnegative integer $n$ is a partition of $n$ where the first occurrence of each distinct part may be overlined.
To define these ranks, let $l(\lambda)$ be the largest part of $\lambda$, $n(\lambda)$ be the number of parts of $\lambda$ and $n(\lambda_o)$  be the number of odd non-overlined parts of $\lambda$.

Analogous to the rank of a partition, the $D$-rank of an overpartition $\lambda$ is the largest part $l(\lambda)$ minus the number of parts $n(\lambda)$.
Similarly,  let $\overline{N}(s,m,n)$ to be the number of overpartitions of $n$ with rank congruent to $s$ modulo $m$. Lovejoy and Osburn \cite{Lovejoy-Osburn-2008}  studied  the rank differences  $\overline{N}(s, \ell, \ell n + d)- \overline{N}(t, \ell, \ell n + d)$ with $\ell = 3$ and $5$ and $0\leq d, s, t < \ell$.

The $M_2$-rank of an overpartition $\lambda$ is defined by
\[M_2\text{-}\mathrm{rank}(\lambda)=\left\lceil\frac{l(\lambda)}{2}\right\rceil-n(\lambda)-n(\lambda_o)-\chi(\lambda),\]
where $\chi(\lambda)=1$ if $l(\lambda)$ is odd and non-overlined and $\chi(\lambda)=0$ otherwise.
Denote by $\overline{N}_2(m,n)$  the number of overpartitions of $n$ with $M_2$-rank equal to $m$. We adopt the notation $\overline{R2}(z;q)$ in \cite{Jennings-2015} to
denote the generating function of $\overline{N}_2(m,n)$, that is,
\begin{align}\label{GF-M2-1}
\overline{R2}(z;q)=\sum_{n=0}^\infty\sum_{m=-\infty}^{\infty}\overline{N}_2(m,n)z^mq^n.
\end{align}
Lovejoy \cite{Lovejoy-2008} obtained the following  generating function:
\begin{eqnarray}\label{GF-M2}
\overline{R2}(z;q)%=\sum_{n=0}^\infty\sum_{m=-\infty}^{\infty}\overline{N}_2(m,n)z^mq^n
=\frac{(-q;q)_\infty}{(q;q)_\infty}\sum_{n=-\infty}^\infty
\frac{(1-z)(1-z^{-1})(-1)^nq^{n^2+2n}}{(1-zq^{2n})(1-z^{-1}q^{2n})}.
\end{eqnarray}

In this paper, we shall study $M_2$-ranks of overpartitions modulo 6 and 10. The main results are presented in Theorems \ref{main-result-6-M2}, \ref{main-result-10-1-M2} and \ref{main-result-10-2-M2} below. We shall adopt the following common notation
\begin{align*}
&\quad (x_1,x_2,\ldots,x_k;q)_\infty=\prod_{n=0}^{\infty}(1-x_1q^n)(1-x_2q^n)\cdots(1-x_kq^n),
\\[5pt]
&\quad j(z;q)=(z;q)_\infty (q/z;q)_\infty (q;q)_\infty,\\[5pt]
&\quad J_{a,m}=j(q^a,q^m),\ \  J_m=(q^m;q^m)_\infty,\ \ \overline{J}_{a,m}=j(-q^a,q^m).
\end{align*}
where we assume that $|q|<1$.

Define
\begin{align*}
\overline{R}_{s,t}(d,\ell)&=\sum_{n=0}^\infty\left(\overline{N}_2(s,\ell,\ell n/2+d)+\overline{N}_2(s+1,\ell,\ell n/2+d)\right.
\\
&\left.\quad\quad\quad\quad-\overline{N}_2(t,\ell,\ell n/2+d)-\overline{N}_2(t+1,\ell,\ell n/2+d)\right)q^n.
\end{align*}

\begin{thm}\label{main-result-6-M2}

We have
\begin{align}\label{main-result-6-M2-1}
\overline{R}_{0,2}(0,6)&=\frac{J_{6}^9J_{2,4}^2}{J^6_{1,6}J^3_{2}J_{3,6}},
\\[3pt] \label{main-result-6-M2-2}
\overline{R}_{0,2}(1,6)&=\frac{2J_{6}^9J_{2,4}^2}{J^5_{1,6}J^3_{2}J^2_{3,6}},
\\[3pt] \label{main-result-6-M2-3}
\overline{R}_{0,2}(2,6)&=\frac{4J_{6}^9J_{2,4}^2}{J^4_{1,6}J^3_{2}J^3_{3,6}}
+\frac{2q^2}{J_{3,6}}
\sum_{n=-\infty}^\infty\frac{(-1)^nq^{3n^2+6n}}{1+q^{6n+4}}.
\end{align}
\end{thm}

\begin{thm}\label{main-result-10-1-M2}
We have
\begin{align}\label{main-result-eq1-10-M2-1}
\overline{R}_{0,4}(i,10)&=\frac{2J_{10}^4J_{5,10}}{J_{1,2}^6}A_i(q),~~0\leq i\leq3,
\\[3pt] \label{main-result-eq1-10-M2-2}
\overline{R}_{0,4}(4,10)&=\frac{2J_{10}^4J_{5,10}}{J_{1,2}^6}A_4(q)
+\frac{q^{4}}{J_{5,10}}
\sum_{n=-\infty}^\infty\frac{(-1)^nq^{5n^2+10n}}{1+q^{10n+8}},
\end{align}
where
\begin{align}\label{10-eq-A0}
A_0(q)&=\frac{J_{2}^4J_{3,20}^2J_{5,20}^8J_{7,20}^2}{2J_{4,20}^4J_{6,20}^3
J_{8,20}J_{20}^6}
+\frac{2qJ_{1}^3J_{10}J_{5,20}J_{20}^5}{J_{1,20}^2J_{2,20}J_{8,20}^3
J_{9,20}^2}
+\frac{8q^{2}J_{1}^2J_{2}^4 J_{2,20}J_{3,20}J_{7,20}}{J_{4,20}^2
J_{6,20}J_{8,20}J_{10,20}^2J_{20}}
\nonumber\\[5pt]&
\quad\quad
-\frac{16q^{2}J_{1}^2 J_{2,20}J_{3,20}^2J_{6,20}J_{7,20}^2}{J_{4,20}
J_{5,20}^2J_{20}^3},
\\[5pt] \label{10-eq-A1}
A_1(q)&=\frac{J_{2}^4 J_{3,20}^3J_{5,20}^6J_{7,20}^3}
{J_{4,20}^4J_{6,20}^3J_{8,20}J_{20}^6}
+\frac{8qJ_{1}^2J_{10}^2J_{2,20}^2J_{3,20}^2J_{7,20}^2J_{8,20}}
{J_{1,20}J_{4,20}^2J_{6,20}J_{9,20}J_{20}^4}
-\frac{8qJ_{1}^2 J_{2,20}^4J_{6,20}^2J_{8,20}^3}
{J_{1,20}^2J_{3,20}J_{7,20}J_{9,20}^2J_{20}^3}
\nonumber\\[5pt]&
\quad\quad
+\frac{4q^{2}J_{1}J_{2}^2 J_{1,20}^2J_{5,20}^3J_{9,20}^2}
{J_{4,20}^2J_{8,20}J_{20}^4J_{10}}
-\frac{8q^{3}J_{1}^2J_{1,20}^2J_{4,20}J_{6,20}J_{9,20}^2J_{10,20}}
{J_{5,20}^2J_{8,20}^2J_{20}^3},
\\[5pt] \label{10-eq-A2}
A_2(q)&=\frac{2J_{2}^4J_{3,20}^4J_{5,20}^4J_{7,20}^4}
{J_{4,20}^4J_{6,20}^3J_{8,20}J_{20}^6}
+\frac{2qJ_{1}^2J_{2,20}J_{5,20}^4J_{6,20}}
{J_{3,20}J_{4,20}J_{7,20}J_{20}^3}
-\frac{8qJ_{1}^2J_{10}^2 J_{3,20}^3J_{4,20}J_{6,20}J_{7,20}^3}
{J_{1,20}J_{5,20}^2J_{8,20}^2J_{9,20}J_{20}^4}
\nonumber\\[5pt]&
\quad\quad
+\frac{4q^{2}J_{1}^2J_{10}^2 J_{1,20}J_{4,20}J_{6,20}J_{9,20}}
{J_{8,20}^2J_{20}^4}
+\frac{8q^{2}J_{1}^2 J_{1,20}J_{2,20}J_{6,20}J_{9,20}}
{J_{4,20}J_{20}^3},
\\[5pt] \label{10-eq-A3}
A_3(q)&=-\frac{J_{1}J_{2}^2J_{5,20}^7}
{J_{10}J_{4,20}^2J_{8,20}J_{20}^4}
+\frac{4J_{2}^4J_{3,20}^5J_{5,20}^2J_{7,20}^5}
{J_{4,20}^4J_{6,20}^3J_{8,20}J_{20}^6}
+\frac{2qJ_{1}^2J_{10}^2 J_{2,20}^2J_{5,20}^2J_{8,20}}
{J_{4,20}^2J_{6,20}J_{20}^4}
\nonumber\\[5pt]&
\quad\quad
-\frac{12qJ_{1}^2 J_{2,20}J_{5,20}^2J_{6,20}}
{J_{4,20}J_{20}^3}
-\frac{16q^{2}J_{1}^3J_{20}^6}{J_{4,20}^2J_{5,20}^3J_{8,20}J_{10}}
+\frac{8q^{2}J_{1}^3J_{10}J_{2,20}J_{20}^5}
{J_{4,20}^3J_{5,20}^3J_{6,20}^2},
\\[5pt] \label{10-eq-A4}
A_4(q)&=-\frac{2J_{1}^2J_{2,20}J_{5,20}^4J_{6,20}}
{J_{1,20}J_{4,20}J_{9,20}J_{20}^3}
-\frac{J_{1}^2J_{10}^2 J_{2,20}^2J_{5,20}^4J_{8,20}}
{J_{1,20}J_{4,20}^2J_{6,20}J_{9,20}J_{20}^4}
+\frac{8J_{2}^4J_{3,20}^6J_{7,20}^6}{J_{4,20}^4J_{6,20}^3J_{8,20}J_{20}^6}
\nonumber\\[5pt]&
\quad\quad
-\frac{16qJ_{1}^2J_{2,20}J_{3,20}J_{6,20}J_{7,20}}
{J_{4,20}J_{20}^3}
-\frac{8qJ_{1}^3J_{10}J_{20}^5J_{2,20}}{J_{1,20}J_{4,20}^3J_{5,20}
J_{6,20}^2J_{9,20}}
-\frac{16q^{3}J_{1}J_{2}^2J_{1,20}^4J_{9,20}^4}
{J_{10}J_{20}^4J_{4,20}^2J_{5,20}J_{8,20}}.
\end{align}
\end{thm}

\begin{thm}\label{main-result-10-2-M2}
We have
\begin{align}\label{main-result-eq2-10-M2-1}
\overline{R}_{1,3}(i,10)&=\frac{2J_{10}^4J_{5,10}}{J_{1,2}^6}B_i(q),~~i=0,2,3,
\\[3pt]\label{main-result-eq2-10-M2-2}
\overline{R}_{1,3}(1,10)&=\frac{2J_{10}^4J_{5,10}}{J_{1,2}^6}B_1(q)
-\frac{q^{4}}{J_{5,10}}
\sum_{n=-\infty}^\infty\frac{(-1)^nq^{5n^2+10n}}{1+q^{10n+6}},
\\[3pt]\label{main-result-eq2-10-M2-3}
\overline{R}_{1,3}(4,10)&=\frac{2J_{10}^4J_{5,10}}{J_{1,2}^6}B_4(q)
-\frac{q^{4}}{J_{5,10}}
\sum_{n=-\infty}^\infty\frac{(-1)^nq^{5n^2+10n}}{1+q^{10n+8}},
\end{align}
where
\begin{align} \label{10-eq-B0}
B_0(q)&=-\frac{qJ_{1}^5J_{1,20}J_{5,20}J_{9,20}}{J_{10}J_{20}J_{4,20}^2J_{6,20}J_{8,20}}
+\frac{4qJ_{1}^4J_{20}^2J_{3,20}J_{7,20}}{J_{1,20}J_{4,20}^2J_{8,20}J_{9,20}J_{10,20}}
-\frac{16q^{2}J_{2}^{\frac{7}{2}}J_{1,2}^{\frac{1}{2}}J_{1,20}J_{2,20}J_{9,20}}
{J_{20}J_{4,20}J_{5,20}^3}
\nonumber\\[5pt]&
\quad\quad
+\frac{16q^{2}J_{1}^2J_{2,20}J_{3,20}^2J_{6,20}J_{7,20}^2}{J_{20}^3 J_{4,20}J_{5,20}^2}
-\frac{4q^{3}J_{1}^4J_{20}^2J_{1,20}J_{9,20}}
{J_{3,20}J_{4,20}J_{7,20}J_{8,20}^2J_{10,20}}
-\frac{8q^{3}J_{1}J_{2}^2J_{1,20}^2J_{5,20}^3J_{9,20}^2}
{J_{10}J_{20}^4J_{4,20}J_{8,20}^2}
\nonumber\\[5pt]&
\quad \quad -\frac{16q^{4}J_{1,20}^2J_{2,20}J_{6,20}J_{9,20}^2J_{1}^2}{J_{20}^3J_{5,20}^2J_{8,20}},
\\[5pt] \label{10-eq-B1}
B_1(q)&=-\frac{2qJ_{1}^2J_{10}^2J_{2,20}J_{5,20}^2}{J_{20}^4J_{4,20}}
+\frac{8qJ_{1}^3J_{2,20}J_{3,20}^3J_{6,20}J_{7,20}^3}{J_{10}J_{4,20}J_{5,20}J_{20}^6}
+\frac{2q^{2}J_{2}^4J_{1,20}J_{2,20}J_{5,20}^2J_{6,20}J_{9,20}}
{J_{20}^6J_{8,20}J_{10,20}}
\nonumber\\[5pt]&
\quad\quad
+\frac{4q^{2}J_{2}^4J_{1,20}^2J_{2,20}J_{6,20}J_{9,20}^2}{J_{20}^6J_{4,20}J_{10,20}}
-\frac{16q^{2}J_{1}^3J_{10}J_{20}^5}{J_{4,20}^2J_{5,20}^3J_{6,20}J_{8,20}}
+\frac{12q^{3}J_{1}^2J_{1,20}^2J_{2,20}J_{6,20}J_{9,20}^2}{J_{20}^3J_{4,20}J_{5,20}^2}
\nonumber\\[5pt]&
\quad\quad
+\frac{14q^{3}J_{1}^2J_{1,20}J_{2,20}J_{6,20}J_{9,20}}{J_{20}^3J_{8,20}}
-\frac{16q^{4}J_{1}^2J_{10}^2J_{1,20}^4J_{2,20}J_{9,20}^4}
{J_{20}^4J_{3,20}J_{4,20}J_{5,20}^4 J_{7,20}},
\\[5pt] \label{10-eq-B2}
B_2(q)&=\frac{qJ_{1}^3J_{2,20}J_{5,20}^5J_{6,20}}{J_{10}J_{20}^6J_{8,20}}
+\frac{4qJ_{1}^2J_{10}^2J_{1,20}^2J_{4,20}^2J_{6,20}^3J_{9,20}^2}
{J_{20}^4J_{2,20}^2J_{3,20}J_{7,20}J_{8,20}^3}
-\frac{8qJ_{2}^4J_{1,20}J_{3,20}^5J_{7,20}^5J_{9,20}}
{J_{20}^6J_{4,20}^4J_{6,20}^3J_{8,20}}
\nonumber\\[5pt]&
\quad\quad
+\frac{12qJ_{1}J_{2}^2J_{3,20}^4J_{7,20}^4}{J_{10}J_{20}^4 J_{4,20}^2J_{5,20}J_{8,20}}
+\frac{4q^{2}J_{1}J_{10}J_{2,20}^6J_{6,20}^2J_{8,20}^5}{J_{5,20}J_{20}^{12}}
-\frac{12q^{2}J_{1}^2J_{2,20}J_{5,20}^2J_{6,20}}{J_{20}^3J_{8,20}}
\nonumber\\[5pt]&
\quad\quad
-\frac{16q^{2}J_{1}^2J_{1,20}J_{2,20}J_{6,20}J_{9,20}}{J_{20}^3J_{4,20}}
-\frac{8q^{3}J_{1}J_{2}^2J_{10}J_{2,20}^4J_{6,20}}{J_{20}^5J_{4,20}J_{5,20}},
\\[5pt] \label{10-eq-B3}
B_3(q)&=\frac{J_{1}^3J_{2,20}J_{5,20}^5J_{6,20}}{J_{10}J_{20}^6J_{4,20}}
-\frac{8qJ_{1}^2J_{10}^2J_{2,20}^3J_{3,20}^2J_{7,20}^2J_{8,20}^2}
{J_{20}^4J_{1,20}J_{4,20}^3J_{6,20}^2J_{9,20}}
+\frac{12qJ_{1}^2J_{2,20}J_{5,20}^2J_{6,20}}{J_{20}^3J_{4,20}}
\nonumber\\[5pt]&
\quad\quad
+\frac{4q^{2}J_{2}^4J_{1,20}^6J_{5,20}^2J_{6,20}J_{9,20}^6}
{J_{20}^6J_{2,20}^4J_{3,20}J_{7,20}J_{8,20}^5}
-\frac{16q^{3}J_{1}^3J_{20}^6J_{3,20}J_{7,20}}
{J_{10}J_{1,20}J_{5,20}^3J_{8,20}^3J_{9,20}}
-\frac{16q^{4}J_{1}J_{2}^2J_{1,20}^4J_{9,20}^4}
{J_{10}J_{20}^4J_{10,20}J_{5,20}J_{8,20}^2},
\\[5pt] \label{10-eq-B4}
B_4(q)&=\frac{2J_{1}^2J_{2,20}J_{5,20}^4J_{6,20}}
{J_{20}^3J_{1,20}J_{4,20}J_{9,20}}
-\frac{2qJ_{1}^2J_{10}^2J_{2,20}^5J_{5,20}^2J_{6,20}^2J_{8,20}}
{J_{20}^4J_{1,20}^3J_{3,20}J_{7,20}J_{9,20}^3}
+\frac{16q^{2}J_{1}^2J_{10}^2J_{1,20}^2J_{2,20}J_{9,20}^2}{J_{20}^4J_{4,20}J_{5,20}^2}
\nonumber\\[5pt]&
\quad\quad
-\frac{16q^{2}J_{1}^2J_{2,20}J_{3,20}^2J_{6,20}J_{7,20}^2}{J_{20}^3J_{5,20}^2J_{8,20}}
+\frac{8q^{3}J_{1}^2J_{1,20}^2J_{2,20}J_{6,20}J_{9,20}^2}
{J_{20}^3J_{3,20}J_{7,20}J_{8,20}}
+\frac{16q^{3}J_{1}J_{2}^2J_{1,20}^4J_{9,20}^4}
{J_{10}J_{20}^4J_{4,20}^2J_{5,20}J_{8,20}}
\nonumber\\[5pt]&
\quad\quad+\frac{32q^{4}J_{1}^2J_{2,20}^3J_{3,20}^2J_{7,20}^2}
{J_{10}^2J_{20}^2J_{4,20}J_{5,20}^2}.
\end{align}
\end{thm}

The generating functions of rank differences often lead to
rank inequalities,  see, for example, Andrews \cite{Andrews-2000}, Garvan \cite{Garvan-1988}, Mao \cite{Mao-2013,Mao-2015}.
In this paper, we obtain the following inequalities between the $M_2$-ranks of overpartitions modulo 6 with the aid of the generating functions in  Theorem \ref{main-result-6-M2}.

\begin{thm} \label{rank-ineq1} For $n\geq 1$, we have
\begin{align}\label{ineq-3-0-1-f}
\overline{N}_2(0,6,3n)&\geq\overline{N}_2(2,6,3n),\\[5pt]
\label{ineq-3-0-3-f}
\overline{N}_2(1,6,3n)&\geq\overline{N}_2(3,6,3n),\\[5pt]
 \label{ineq-3-1-1-f}
\overline{N}_2(0,6,3n+1)&\geq\overline{N}_2(2,6,3n+1),
\\[5pt] \label{ineq-3-1-3-f}
\overline{N}_2(1,6,3n+1)&\geq\overline{N}_2(3,6,3n+1),
\\[5pt] \label{ineq-3-0-2-f}
\overline{N}_2(0,6,3n)+\overline{N}_2(3,6,3n)&\geq\overline{N}_2(1,6,3n)+\overline{N}_2(2,6,3n),\\[5pt]  \label{ineq-3-1-2-f}
\overline{N}_2(0,6,3n+1)+\overline{N}_2(3,6,3n+1)&\geq\overline{N}_2(1,6,3n+1)+\overline{N}_2(2,6,3n+1),
\\[5pt] \label{ineq-3-2-1-f}
\overline{N}_2(0,6,3n+2)+\overline{N}_2(3,6,3n+2)&\geq\overline{N}_2(1,6,3n+2)+\overline{N}_2(2,6,3n+2).
\end{align}
\end{thm}

The connection between classical mock theta functions and the generating functions of rank differences of partitions have been extensively studied, see \cite{Andrews-Garvan-1989,Hickerson-1988,Lovejoy-Osburn-2010}. We find a relation between the generating function of $M_2$-ranks of overpartitions modulo 6 and the third order mock theta function $\rho(q)$ defined by
$$
\rho(q)=\sum_{n=0}^\infty\frac{q^{2n(n+1)}(q;q^2)_{n+1}}{(q^3;q^6)_{n+1}},$$ see \cite{Watson-1936}.

\begin{thm}\label{mock-1}We have
\begin{align}\label{mock-1-eqn}
\overline{R}_{0,2}(2,6)=2\rho(-q)+M_0(q),
\end{align}
where $M_0(q)$ is an (explicit) weakly holomorphic modular form given by
\[M_0(q)=-\frac{2 J^3_{6}J_2\overline{J}_{1,6}}{J_{1,6}J_{3,6}\overline{J}_{2}^2}
+\frac{4J_{6}^9J_{2,4}^2}{J^4_{1,6}J^3_{2}J^3_{3,6}}.\]
\end{thm}

In light of Theorem \ref{main-result-10-1-M2} and Theorem \ref{main-result-10-2-M2}, we deduce the following relations on the generating function of $M_2$-ranks of overpartitions modulo 10 and the tenth order mock theta functions $\phi(q)$ and $\psi(q)$ are defined as follows,
\begin{align*}
  \phi(q)=\sum_{n=0}^\infty\frac{q^{n+1\choose 2}}{(q;q^2)_{n+1}},  \quad \quad
  \psi(q)=\sum_{n=0}^\infty\frac{q^{n+2\choose 2}}{(q;q^2)_{n+1}},
\end{align*}
see \cite{Choi-1999}.

\begin{thm}\label{mock-2}We have
\begin{align}
\label{mock-2-1}
\overline{R}_{0,4}(4,10)&=q^{-1}\psi(-q)+M_1(q),
\\[5pt] \label{mock-2-2}
\overline{R}_{1,3}(1,10)
&=\phi(-q)+M_2(q),
\\[5pt] \label{mock-2-3}
\overline{R}_{1,3}(4,10)
&=-q^{-1}\psi(-q)+M_3(q),
\end{align}
where $M_1(q),\ M_2(q)$ and $M_3(q)$ are (explicit) weakly holomorphic modular forms given by
\begin{align*}
M_1(q)&=\frac{2J_{10}^4J_{5,10}}{J_{1,2}^6}A_4(q)
    +\frac{J_{10}^3 J_{1,10}\overline{J}_{4,10}}
{J_{4,10}J_{5,10}\overline{J}_{1,10}\overline{J}_{2,10}}
  -\frac{qJ_{10}^3J_{1,10}\overline{J}_{2,10}}
  {J_{4,10}J_{5,10}\overline{J}_{2,10}\overline{J}_{3,10}},
\\[3pt]
M_2(q)&=\frac{2J_{10}^4J_{5,10}}{J_{1,2}^6}B_1(q)
    -\frac{J_{10}^3J_{3,10}\overline{J}_{4,10}}{J_{2,10}J_{5,10}\overline{J}_{1,10}\overline{J}_{4,10}}
    +\frac{qJ_{10}^3J_{3,10}\overline{J}_{2,10}}
    {J_{2,10}J_{5,10}\overline{J}_{3,10}\overline{J}_{4,10}}.
\\[3pt]
M_3(q)&= \frac{2J_{10}^4J_{5,10}}{J_{1,2}^6}B_4(q)
    -\frac{J_{10}^3\overline{J}_{4,10}J_{1,10}}
{J_{4,10}J_{5,10}\overline{J}_{1,10}\overline{J}_{2,10}}
  +\frac{qJ_{10}^3J_{1,10}\overline{J}_{2,10}}
  {J_{4,10}J_{5,10}\overline{J}_{2,10}\overline{J}_{3,10}},
\end{align*}
and $A_4(q)$ is defined in Theorem \ref{main-result-10-1-M2}, $B_1(q)$ and $B_4(q)$ are defined in Theorem \ref{main-result-10-2-M2}.
\end{thm}

This paper is organized as follows. We prove Theorem \ref{main-result-6-M2} in Section 2.
Theorems \ref{main-result-10-1-M2} and  \ref{main-result-10-2-M2} are proved in Section 3.
In Section 4, we derive inequalities on $M_2$-ranks of overpartitions mod 6 from the generating function in Theorem \ref{main-result-6-M2}. In Section 5, we obtain relations between the generating functions of $M_2$-rank differences of overpartitions and mock theta functions as given in Theorems \ref{mock-1} and \ref{mock-2}.

\section{Proof of Theorem \ref{main-result-6-M2}}

The following lemma plays a central role in the proof of Theorem \ref{main-result-6-M2}.

\begin{lem}\label{main-result-6-lem-M2}
We have
\begin{align}
\overline{R2}(\exp(\pi i/3);q)&=\sum_{n=0}^\infty(\overline{N}_2(0,6,n)+\overline{N}_2(1,6,n)
-\overline{N}_2(2,6,n)-\overline{N}_2(3,6,n))q^n
\nonumber\\
&=\frac{2(-q;q)_\infty}{(q;q)_\infty}\sum_{n=-\infty}^\infty
\frac{(-1)^nq^{n^2+2n}}{1+q^{6n}}.
\end{align}
\end{lem}

\pf Replacing $z$ by $\xi_6=e^{\frac{\pi i}{3}}$, the right-hand side of \eqref{GF-M2-1} simplifies to
\begin{align}\label{R2---1}
\overline{R2}(\exp(\pi i/3);q)&=\sum_{n=0}^\infty\sum_{m=-\infty}^{\infty}\overline{N}_2(m,n)\xi_6^mq^n
\nonumber\\
&=\sum_{n=0}^\infty\sum_{t=0}^5\sum_{m=-\infty}^{\infty}\overline{N}_2(6m+t,n)\xi_6^tq^n
\nonumber\\
&=\sum_{n=0}^\infty\left(\overline{N}_2(0,6,n)+\xi_6\overline{N}_2(1,6,n)
+\xi_6^2\overline{N}_2(2,6,n)\right.
\nonumber\\&
\ \ \ \ \ \ \
\left.+\xi_6^3\overline{N}_2(3,6,n)+\xi^4_6\overline{N}_2(4,6,n)+\xi_6^5\overline{N}_2(5,6,n)\right)q^n.
\end{align}
In view of the symmetry property
\[\overline{N}_2(m,n)=\overline{N}_2(-m,n),~\text{for}~m\geq0\]
due to Lovejoy and Osburn \cite{Lovejoy-Osburn-2010}, we see that
\begin{align}\label{Sy-N}
\overline{N}_2(s,m,n)=\overline{N}_2(m-s,m,n).
\end{align}
In particular, for $m\geq0$,
\begin{align}\label{N(s,m,n)-1}
\overline{N}_2(1,6,n)=\overline{N}_2(5,6,n),
\end{align}
and
\begin{align}\label{N(s,m,n)-2}
\overline{N}_2(2,6,n)=\overline{N}_2(4,6,n).
\end{align}
Substituting \eqref{N(s,m,n)-1} and \eqref{N(s,m,n)-2} into \eqref{R2---1}, we find that
\begin{align}\label{R2----2}
\overline{R2}(\exp( \pi i/3);q)&=\sum_{n=0}^\infty\left(\overline{N}_2(0,6,n)+(\xi_6+\xi_6^5)\overline{N}_2(1,6,n)\right.
\nonumber\\
&\quad\quad\quad\left.+(\xi_6^2+\xi^4_6)\overline{N}_2(2,6,n)+\xi_6^3\overline{N}_2(3,6,n)\right)q^n.
\end{align}
Observing that $1-\xi_6+\xi_6^2=0$ and $\xi_6^3=-1$, \eqref{R2----2} yields that
\begin{align*}
\overline{R2}(\exp( \pi i/3);q)=\sum_{n=0}^\infty(\overline{N}_2(0,6,n)+\overline{N}_2(1,6,n)
-\overline{N}_2(2,6,n)-\overline{N}_2(3,6,n))q^n.
\end{align*}
To evaluate the right-hand side of \eqref{GF-M2} for $z=\xi_6$, we get
\begin{align*}
\overline{R2}(\exp( \pi i/3);q)
&=\frac{(-q;q)_\infty}{(q;q)_\infty}\sum_{n=-\infty}^\infty
\frac{(1-\xi_6)(1-\xi_6^{-1})(-1)^nq^{n^2+2n}}{(1-\xi_6q^{2n})(1-\xi_6^{-1}q^{2n})}
\\
&=\frac{(-q;q)_\infty}{(q;q)_\infty}\sum_{n=-\infty}^\infty
\frac{(2-\xi_6^{-1}-\xi_6)(-1)^nq^{n^2+2n}}{1-\xi_6^{-1}q^{2n}-\xi_6q^{2n}+q^{4n}}.
\end{align*}
But $1-\xi_6^{-1}-\xi_6=0$, so we obtain that
\begin{align*}
\overline{R2}(\exp( \pi i/3);q)&=\frac{(-q;q)_\infty}{(q;q)_\infty}\sum_{n=-\infty}^\infty
\frac{(-1)^nq^{n^2+2n}}
{1-q^{2n}+q^{4n}}
\\
&=\frac{(-q;q)_\infty}{(q;q)_\infty}\sum_{n=-\infty}^\infty
\frac{(-1)^nq^{n^2+2n}(1+q^{2n})}{1+q^{6n}}
\\
&=\frac{2(-q;q)_\infty}{(q;q)_\infty}\sum_{n=-\infty}^\infty
\frac{(-1)^nq^{n^2+2n}}{1+q^{6n}},
\end{align*}
as desired. This completes the proof of Lemma \ref{main-result-6-lem-M2}. \qed

Now we are in a position to prove Theorem  \ref{main-result-6-M2}. \vskip 0.2cm

\noindent{\it Proof of Theorem  \ref{main-result-6-M2}.} By Lemma \ref{main-result-6-lem-M2}, it suffices to show that
\begin{align}\label{mod6-eq1-M2}
&\frac{(-q;q)_\infty}{(q;q)_\infty}\sum_{n=-\infty}^\infty
\frac{(-1)^nq^{n^2+2n}}{1+q^{6n}}\nonumber \\[5pt]
&\quad=\frac{J_{18}^9J_{6,12}^2}{2J^6_{3,18}J^3_{6}J_{9,18}}
+q\frac{J_{18}^9J_{6,12}^2}{J^5_{3,18}J^3_{6}J^2_{9,18}}
\nonumber\\&\quad\quad\quad 
+q^2\left\{\frac{2J_{18}^9J_{6,12}^2}{J^4_{3,18}J^3_{6}J^3_{9,18}}
+\frac{q^6}{J_{9,18}}
\sum_{n=-\infty}^\infty\frac{(-1)^nq^{9n^2+18n}}{1+q^{18n+12}}\right\}.
\end{align}
The summation on the left hand side of \eqref{mod6-eq1-M2} can be written as
\begin{align}\label{main-result-6-disse-M2}
\sum_{n=-\infty}^\infty\frac{(-1)^nq^{n^2+2n}}{1+q^{6n}}
&=\sum_{n=-\infty}^\infty\frac{(-1)^nq^{9n^2+6n}}{1+q^{18n}}
-\sum_{n=-\infty}^\infty\frac{(-1)^nq^{9n^2+12n+3}}{1+q^{18n+6}}
+\sum_{n=-\infty}^\infty\frac{(-1)^nq^{9n^2+18n+8}}{1+q^{18n+12}}.
\end{align}
Let $G_0, G_1$ and $G_2$ denote the summations on the right hand side
of \eqref{main-result-6-disse-M2}, respectively.
Recall the following  identity   \cite[Lemma 3.1]{Lovejoy-Osburn-2010}:
\begin{align}\label{Lovejoy-Osburn-2010}
&\sum_{n=-\infty}^\infty(-1)^nq^{n^2+2n}\left[\frac{\zeta^{-2n}}{1-z^2\zeta^{-2}q^{2n}}
+\frac{\zeta^{2n+4}}{1-z^2\zeta^2q^{2n}}\right]
\nonumber\\
&\quad=\frac{-2(\zeta^4,q^2\zeta^{-4};q^2)_\infty(-q;q)_\infty^2}
{(-\zeta^2,-q\zeta^{-2};q)_\infty(\zeta^{-2},q^2\zeta^2;q^2)_\infty}
\sum_{n=-\infty}^\infty(-1)^n\frac{q^{n^2+2n}}{1-z^2q^{2n}}
\nonumber\\
&\ \ \ \
\quad+\frac{(-z^2,-qz^{-2};q)_\infty
(\zeta^4,q^2\zeta^{-4},\zeta^2,q^2\zeta^{-2};q^2)_\infty(q^2;q^2)_\infty^2}
{(-\zeta^2,-q\zeta^{-2};q)_\infty
(z^2\zeta^2,q^2z^{-2}\zeta^{-2},z^2\zeta^{-2},q^2z^{-2}\zeta^2,z^2,q^2z^{-2};q^2)_\infty}.
\end{align}
Substituting $q$, $z^2$ and $\zeta^2$  by $q^9$, $-q^{12}$ and $q^{12}$, respectively, in \eqref{Lovejoy-Osburn-2010}, we find that
\begin{align*}
&\sum_{n=-\infty}^\infty(-1)^nq^{9n^2+18n}
\left[\frac{q^{-12n}}{1+q^{18n}}+\frac{q^{12n+24}}{1+q^{18n+24}}\right]
\\
&\quad=\sum_{n=-\infty}^\infty
\frac{(-1)^nq^{9n^2+6n}}{1+q^{18n}}-\sum_{n=-\infty}^\infty\frac{(-1)^nq^{9n^2+12n+3}}{1+q^{18n+6}}
\\
&\quad=-\frac{2q^9J_{3,18}}{J_{9,18}}
\sum_{n=-\infty}^\infty\frac{(-1)^nq^{9n^2+18n}}{1+q^{18n+12}}
%\\&
+\frac{J_{6}J_{3,18}^2J_{6,12}^2}{2J_{18}^3},
\end{align*}
which can be recast as
\begin{equation}\label{main-result-6-claim-M2}
G_0-G_1=-\frac{2qJ_{3,18}}{J_{9,18}}G_2
+\frac{J_{6}J_{3,18}^2J_{6,12}^2}{2J_{18}^3}.
\end{equation}
Combining \eqref{main-result-6-claim-M2} into \eqref{main-result-6-disse-M2}, we are led to the relation
\begin{align}\label{main-result-6-claim-M2--1}
\sum_{n=-\infty}^\infty\frac{(-1)^nq^{n^2+2n}}{1+q^{6n}}
%&=F_0-F_1+F_2\\
&=\left(1-\frac{2qJ_{3,18}}{J_{9,18}}\right)G_2
+\frac{J_{6}J_{3,18}^2J_{6,12}^2}{2J_{18}^3}.
\end{align}
Applying the following identity \cite[Lemma 3.1]{Lovejoy-Osburn-2008}
\begin{align}\label{Lovejoy-Osburn-2008-mod6}
\frac{(q;q)_\infty}{(-q;q)_\infty}=\frac{(q^9;q^9)_\infty}{(-q^9;q^9)_\infty}
-2q(q^3,q^{15},q^{18};q^{18})_\infty
=J_{9,18}-2qJ_{3,18}
\end{align}
to \eqref{main-result-6-claim-M2--1}, we get
\begin{align*}
\sum_{n=-\infty}^\infty\frac{(-1)^nq^{n^2+2n}}{1+q^{6n}}
&=\frac{(q;q)_\infty}{(-q;q)_\infty}\cdot\frac{G_2}{J_{9,18}}
+\frac{J_{6}J_{3,18}^2J_{6,12}^2}{2J_{18}^3}.
\end{align*}
To complete the proof of Theorem \ref{main-result-6-M2}, it remains to show that
\begin{align}\label{V-eq-M2}
\frac{J_{6}J_{3,18}^2J_{6,12}^2}{J_{18}^3}=\frac{(q;q)_\infty}{(-q;q)_\infty}
\left\{\frac{J_{18}^9J_{6,12}^2}{J^3_{6}J^6_{3,18}J_{9,18}}
+q\frac{2J_{18}^9J_{6,12}^2}{J^3_{6}J^5_{3,18}J^2_{9,18}}
+q^2\frac{4J_{18}^9J_{6,12}^2}{J^3_{6}J^4_{3,18}J^3_{9,18}}\right\}.
\end{align}
 In view of \eqref{Lovejoy-Osburn-2008-mod6}, the right hand side of \eqref{V-eq-M2}
can be rewritten as
\begin{align*}
\left(J_{9,18}-2qJ_{3,18}\right)\left\{\frac{J_{18}^9J_{6,12}^2}{J^3_{6}J^6_{3,18}J_{9,18}}
+q\frac{2J_{18}^9J_{6,12}^2}{J^3_{6}J^5_{3,18}J^2_{9,18}}
+q^2\frac{4J_{18}^9J_{6,12}^2}{J^3_{6}J^4_{3,18}J^3_{9,18}}\right\},
\end{align*}
which equals
\begin{align*}
&\left(\frac{J_{18}^9J_{6,12}^2}{J^3_{6}J^6_{3,18}}
-8q^3\frac{J_{18}^9J_{6,12}^2}{J^3_{6}J^3_{3,18}J^3_{9,18}}\right)
+q\left(\frac{2J_{18}^9J_{6,12}^2}{J^3_{6}J^5_{3,18}J_{9,18}}
-\frac{2J_{18}^9J_{6,12}^2}{J^3_{6}J^5_{3,18}J_{9,18}}\right)
\\[5pt]
&\quad +q^2\left(\frac{4J_{18}^9J_{6,12}^2}{J^3_{6}J^4_{3,18}J^2_{9,18}}
-\frac{4J_{18}^9J_{6,12}^2}{J^3_{6}J^4_{3,18}J^2_{9,18}}\right)\nonumber \\[5pt]
&=\frac{J_{18}^9J_{6,12}^2}{J^3_{6}J^6_{3,18}}
-8q^3\frac{J_{18}^9J_{6,12}^2}{J^3_{6}J^3_{3,18}J^3_{9,18}}.
\end{align*}
Hence, \eqref{V-eq-M2} is equivalent to
\begin{align}\label{U-0}
\frac{J_{18}^9J_{6,12}^2}{J^3_{6}J^6_{3,18}}
-8q^3\frac{J_{18}^9J_{6,12}^2}{J^3_{6}J^3_{3,18}J^3_{9,18}}=\frac{J_{6}J_{3,18}^2J_{6,12}^2}{J_{18}^3},
\end{align}
which can be simplified to
\begin{equation}\label{V-0-r}
J^3_{9,18}-8q^3J^3_{3,18}=\frac{J_{6}^4 J_{3,18}^8J_{9,18}^3}{J_{18}^{12}}.
\end{equation}
The above identity can be deduced from the following identity \cite{Atkin-Swinnerton-Dyer-1954} by setting  $q=q^9$, $x=-q^3$, $y=q^3$ and $z=-1$:
\begin{equation*}%\label{jac}
j(x,q)^2j(yz,q)j(yz^{-1},q)=j(y,q)^2j(xz,q)j(xz^{-1},q)-yz^{-1}j(z,q)^2j(xy,q)j(x y^{-1},q).
\end{equation*}
This completes the proof of Theorem \ref{main-result-6-M2}. \qed

\section{Proofs of Theorem \ref{main-result-10-1-M2} and Theorem \ref{main-result-10-2-M2}}

To prove Theorems \ref{main-result-10-1-M2} and  \ref{main-result-10-2-M2}, we need the following expression for  $\overline{R2}(\exp(\pi i/5);q)$.

\begin{lem}\label{main-result-10-pf-lem-M2}We have
\begin{align}\label{main-result-10-lem-eq-M2}
&\overline{R2}(\exp(\pi i/5);q)
=F_1(q)+(\xi_{10}^2-\xi_{10}^3)F_2(q),
\end{align}
where
\begin{align}\label{M2-10-1}
&F_1(q)=\frac{2(-q;q)_\infty}{(q;q)_\infty}\sum_{n=-\infty}^\infty
\frac{(-1)^nq^{n^2+2n}}{1+q^{10n}},
\\ \label{M2-10-2}
&F_2(q)=\frac{2(-q;q)_\infty}{(q;q)_\infty}\sum_{n=-\infty}^\infty
\frac{(-1)^nq^{n^2+2n}(q^{2n}-1)}{1+q^{10n}}.
\end{align}
\end{lem}

\pf
With the substitution $z=\xi_{10}$,  the right-hand side of \eqref{GF-M2}
equals
\begin{align}\label{10-eq-M2-1}
&\frac{(-q;q)_\infty}{(q;q)_\infty}\sum_{n=-\infty}^\infty
\frac{(1-\xi_{10})(1-\xi_{10}^{-1})(-1)^nq^{n^2+2n}}{(1-\xi_{10}q^{2n})(1-\xi_{10}^{-1}q^{2n})}
\nonumber\\[3pt]
&\quad=\frac{(-q;q)_\infty}{(q;q)_\infty}\sum_{n=-\infty}^\infty
\frac{(1-\xi_{10})(1-\xi_{10}^{-1})(1-\xi_{10}^{3}q^{2n})(1-\xi_{10}^{-3}q^{2n})(1+q^{2n})(-1)^nq^{n^2+2n}}
{(1-\xi_{10}q^{2n})(1-\xi_{10}^{-1}q^{2n})(1-\xi_{10}^{3}q^{2n})(1-\xi_{10}^{-3}q^{2n})(1+q^{2n})}.
\end{align}
Since $1-\xi_{10}-\xi_{10}^{-1}-\xi_{10}^3-\xi_{10}^{-3}=0$,   we deduce that \[(1-\xi_{10}q^{2n})(1-\xi_{10}^{-1}q^{2n})(1-\xi_{10}^{3}q^{2n})(1-\xi_{10}^{-3}q^{2n})(1+q^{2n})=1+q^{10n}\]
and
\begin{align*}
&(1-\xi_{10})(1-\xi_{10}^{-1})(1-\xi_{10}^{3}q^{2n})(1-\xi_{10}^{-3}q^{2n})(1+q^{2n})
\\[3pt]
&=(1-\xi_{10}^2+\xi_{10}^3)+(-1+\xi_{10}^{-1}+\xi_{10})q^{2n}+(-1+\xi_{10}^{-1}+\xi_{10})q^{4n}
+(1-\xi_{10}^2+\xi_{10}^3)q^{6n}
\\[3pt]
&=(1-\xi_{10}^2+\xi_{10}^3)+(\xi_{10}^2-\xi_{10}^3)q^{2n}+(\xi_{10}^2-\xi_{10}^3)q^{4n}+(1-\xi_{10}^2+\xi_{10}^3)q^{6n}.
\end{align*}
Thus the right hand side of \eqref{10-eq-M2-1} equals
\begin{align*}
&\overline{R2}(\exp(\pi i/5);q)\\[3pt]
&=\frac{(-q;q)_\infty}{(q;q)_\infty}\left\{\sum_{n=-\infty}^\infty
\frac{(-1)^nq^{n^2+2n}(1+q^{6n})}
{1+q^{10n}}\right.
\\[3pt]
&\quad \quad \quad \quad \quad \quad \left.+(\xi_{10}^2-\xi_{10}^3)\sum_{n=-\infty}^\infty
\frac{(-1)^nq^{n^2+2n}(q^{2n}+q^{4n}-1-q^{6n})}{1+q^{10n}}\right\}
\\[3pt]
&=\frac{(-q;q)_\infty}{(q;q)_\infty}\left\{\sum_{n=-\infty}^\infty
\frac{(-1)^nq^{n^2+2n}}{1+q^{10n}}+\sum_{n=-\infty}^\infty
\frac{(-1)^{-n}q^{n^2-8n}}{1+q^{-10n}}\right\}
\\[3pt]
&\quad\quad\quad +(\xi_{10}^2-\xi_{10}^3)\frac{(-q;q)_\infty}{(q;q)_\infty}
\left\{\sum_{n=-\infty}^\infty\frac{(-1)^nq^{n^2+4n}}{1+q^{10n}}
+\sum_{n=-\infty}^\infty\frac{(-1)^{-n}q^{n^2-6n}}{1+q^{-10n}}\right.
\\[3pt]
&\quad \quad \quad \quad \quad \quad \quad \quad \quad \quad \quad
\left.-\sum_{n=-\infty}^\infty
\frac{(-1)^nq^{n^2+2n}}{1+q^{10n}}-\sum_{n=-\infty}^\infty
\frac{(-1)^{-n}q^{n^2-8n}}{1+q^{-10n}}\right\}
\\[3pt]
&=\frac{2(-q;q)_\infty}{(q;q)_\infty}\sum_{n=-\infty}^\infty
\frac{(-1)^nq^{n^2+2n}}{1+q^{10n}}
\\[3pt]
&\quad -(\xi_{10}^2-\xi_{10}^3)\frac{2(-q;q)_\infty}{(q;q)_\infty}\sum_{n=-\infty}^\infty
\frac{(-1)^nq^{n^2+2n}(1-q^{2n})}{1+q^{10n}}.
\end{align*}
This completes the proof. \qed

\begin{lem}\label{main-result-10-lem-1}
We have
\begin{align}\label{main-result-10-M2-F1}
&\sum_{n=0}^\infty(\overline{N}_2(0,10,n)+\overline{N}_2(1,10,n)
-\overline{N}_2(4,10,n)-\overline{N}_2(5,10,n))q^n
=F_1(q)
\end{align}
and
\begin{align}\label{main-result-10-M2-F2}
&\sum_{n=0}^\infty(\overline{N}_2(1,10,n)+\overline{N}_2(2,10,n)
-\overline{N}_2(3,10,n)-\overline{N}_2(4,10,n))q^n
=F_2(q),
\end{align}
where $F_1(q)$ and $F_2(q)$ are given in \eqref{M2-10-1} and \eqref{M2-10-2}.
\end{lem}

\pf Setting $z=\xi_{10}=e^{\frac{\pi i}{5}}$, we get
\begin{align*}
\overline{R2}(\exp(\pi i/5);q)&=\sum_{n=0}^\infty\sum_{m=-\infty}^{\infty}\overline{N}_2(m,n)\xi_{10}^mq^n
\\[3pt]
&=\sum_{n=0}^\infty\sum_{t=0}^{9}\sum_{m=-\infty}^{\infty}\overline{N}_2(10m+t,n)
\xi_{10}^tq^n
\\[3pt]
&=\sum_{n=0}^\infty\sum_{t=0}^{9}\overline{N}_2(t,10,n)
\xi_{10}^tq^n.
\end{align*}
Utilizing \eqref{Sy-N} and the fact that $\xi_{10}^5=-1$, we deduce that
\begin{align*}
&\overline{R2}(\exp(\pi i/5);q)
\\[5pt]
&\quad=\sum_{n=0}^\infty\left(\overline{N}_2(0,10,n)+(\xi_{10}-\xi_{10}^4)\overline{N}_2(1,10,n)
+(\xi_{10}^2-\xi_{10}^3)\overline{N}_2(2,10,n)\right.
\\[5pt]
&\quad \quad\quad\quad \left.+(\xi_{10}^3-\xi_{10}^2)\overline{N}_2(3,10,n)+(\xi_{10}^4-\xi_{10})\overline{N}_2(4,10,n)-\overline{N}_2(5,10,n)\right)
q^n.
\end{align*}
But $1-\xi_{10}+\xi_{10}^2-\xi_{10}^3+\xi_{10}^4=0$, so we arrive at
\begin{align}\label{R2-10--1}
&\overline{R2}(\exp(\pi i/5);q)
\nonumber\\[3pt]
&\quad=\sum_{n=0}^\infty\left(\overline{N}_2(0,10,n)+(1+\xi_{10}^2-\xi_{10}^3)\overline{N}_2(1,10,n)
+(\xi_{10}^2-\xi_{10}^3)\overline{N}_2(2,10,n)\right.
\nonumber\\[3pt]
& \quad\quad \quad \quad \left.+(\xi_{10}^3-\xi_{10}^2)\overline{N}_2(3,10,n)
-(1+\xi_{10}^2-\xi_{10}^3)\overline{N}_2(4,10,n)-\overline{N}_2(5,10,n)\right)
q^n
\nonumber\\[5pt]
&\quad=\sum_{n=0}^\infty(\overline{N}_2(0,10,n)+\overline{N}_2(1,10,n)
-\overline{N}_2(4,10,n)-\overline{N}_2(5,10,n))q^n
\nonumber\\
&\quad \quad \quad \quad +(\xi_{10}^2-\xi_{10}^3)\sum_{n=0}^\infty(\overline{N}_2(1,10,n)+\overline{N}_2(2,10,n)
-\overline{N}_2(3,10,n)-\overline{N}_2(4,10,n))q^n.
\end{align}
Combining \eqref{main-result-10-lem-eq-M2} with \eqref{R2-10--1}, we have
\begin{align}\label{10-xi}
&F_1(q)-\sum_{n=0}^\infty(\overline{N}_2(0,10,n)+\overline{N}_2(1,10,n)
-\overline{N}_2(4,10,n)-\overline{N}_2(5,10,n))q^n
\nonumber\\[3pt]
&=(\xi_{10}^2-\xi_{10}^3)\left(\sum_{n=0}^\infty(\overline{N}_2(1,10,n)+\overline{N}_2(2,10,n)
-\overline{N}_2(3,10,n)-\overline{N}_2(4,10,n))q^n-F_2(q)\right).
\end{align}
Let us rewrite \eqref{10-xi} as
\begin{align}\label{10-xi-1}
A(q)=(\xi_{10}^2-\xi_{10}^3)B(q).
\end{align}
where
\[A(q)=\sum_{n=0}^\infty a_nq^n = F_1(q)-\sum_{n=0}^\infty(\overline{N}_2(0,10,n)+\overline{N}_2(1,10,n)
-\overline{N}_2(4,10,n)-\overline{N}_2(5,10,n))q^n,\]
and
\[B(q)=\sum_{n=0}^\infty b_nq^n
=\sum_{n=0}^\infty(\overline{N}_2(1,10,n)+\overline{N}_2(2,10,n)
-\overline{N}_2(3,10,n)-\overline{N}_2(4,10,n))q^n-F_2(q).\]
From \eqref{10-xi-1} it follows that for $n\geq 0$,
\[a_n=(\xi_{10}^2-\xi_{10}^3)b_n.\]
It is known that the minimal polynomial over $\mathbb{Q}$ of $\xi_{10}$ is
$1+\xi_{10}+\xi_{10}^2+\xi_{10}^3+\xi_{10}^4$.
This implies that $\xi_{10}^2-\xi_{10}^3$ is not a rational number.
On the other hand, $a_n$ and $b_n$ are integers for $n\geq 0$.
So we conclude that $a_n=b_n=0$ for $n\geq 0$.
This completes the proof.
\qed

The following lemmas will be used
 in the proofs of  Theorems \ref{main-result-10-1-M2} and \ref{main-result-10-2-M2}.

\begin{lem}\label{mod10-lem-M2}
Let
\begin{align*}
&U_1=\frac{J_{10,100}^4J_{15,100}^2J_{30,100}J_{35,100}^2J_{40,100}^3J_{50,100}^2}
{J_{100}^{12}},
\\
&U_2=\frac{J_{5,100}J_{10,100}^3J_{15,100}J_{20,100}J_{30,100}^3
J_{35,100}J_{40,100}^2J_{45,100}J_{50,100}}{J_{100}^{12}}.
\end{align*}
We have
\begin{align}\label{mod10-M2-eq1}
\sum_{n=-\infty}^\infty
\frac{(-1)^nq^{n^2+2n}}{1+q^{10n}}=\frac{1}{2}U_1-q^3U_2+
\frac{(q;q)_\infty}{(-q;q)_\infty}\cdot\frac{q^{24}}{J_{25,50}}
\sum_{n=-\infty}^\infty\frac{(-1)^nq^{25n^2+50n}}{1+q^{50n+40}}.
\end{align}

\end{lem}

\pf
First, we split the sum of the left side of \eqref{mod10-M2-eq1} as five sums:
\begin{align}\label{main-result-10-disse-1-M2}
\sum_{n=-\infty}^\infty\frac{(-1)^nq^{n^2+2n}}{1+q^{10n}}
&=\sum_{n=-\infty}^\infty\frac{(-1)^nq^{25n^2+10n}}{1+q^{50n}}
-\sum_{n=-\infty}^\infty\frac{(-1)^nq^{25n^2+20n+3}}{1+q^{50n+10}}
+\sum_{n=-\infty}^\infty\frac{(-1)^nq^{25n^2+30n+8}}{1+q^{50n+20}}
\nonumber\\
&\quad\quad\quad -\sum_{n=-\infty}^\infty\frac{(-1)^nq^{25n^2+40n+15}}{1+q^{50n+30}}
+\sum_{n=-\infty}^\infty\frac{(-1)^nq^{25n^2+50n+24}}{1+q^{50n+40}}.
\end{align}
Let $P_0, P_1, P_2, P_3$ and $P_4$ denote the summations on the right hand side
of \eqref{main-result-10-disse-1-M2}, respectively.
 Substituting $q$, $z^2$ and $\zeta^2$  by $q^{25}$, $-q^{40}$ and $q^{40}$, respectively, \eqref{Lovejoy-Osburn-2010} becomes
\begin{align}\label{10-eq---1}
&\sum_{n=-\infty}^\infty\frac{(-1)^nq^{25n^2+10n}}{1+q^{50n}}
-\sum_{n=-\infty}^\infty\frac{(-1)^nq^{25n^2+90n+80}}{1+q^{50n+80}}
\nonumber\\&
\quad\quad=-\frac{2q(q^{20},q^{30};q^{50})_\infty(-q^{25},q^{25})_\infty^2}
{(-q^{10},-q^{15};q^{25})_\infty(q^{10},q^{40};q^{50})_\infty}
\sum_{n=-\infty}^\infty\frac{(-1)^nq^{25n^2+50n+24}}{1+q^{50n+40}}
\nonumber\\
&\quad\quad\quad\quad
+\frac{(q^{10},q^{15};q^{25})_\infty(q^{10};q^{10})_\infty(q^{50};q^{50})_\infty}
{2(-q^{10},-q^{15};q^{25})_\infty(-q^{10};q^{10})_\infty(-q^{50};q^{50})_\infty}.
\end{align}
Replacing $n$ by $n-1$ in second summation on the left side of \eqref{10-eq---1}, we get
\begin{align*}
&\sum_{n=-\infty}^\infty\frac{(-1)^nq^{25n^2+10n}}{1+q^{50n}}
-\sum_{n=-\infty}^\infty\frac{(-1)^nq^{25n^2+40n+15}}{1+q^{50n+30}}
\\
&\quad=-\frac{2qJ_{15,50}}{J_{25,50}}\sum_{n=-\infty}^\infty\frac{(-1)^nq^{25n^2+50n+24}}{1+q^{50n+40}}
+\frac{1}{2}U_1,
\end{align*}
which can be expressed as
\begin{align}\label{main-result-10-claim-1-M2}
&P_0-P_3=-\frac{2qJ_{15,50}}{J_{25,50}}P_4+\frac{1}{2}U_1.
\end{align}
Similarly, substituting  $q$, $z^2$ and $\zeta^2$  by $q^{25}$, $-q^{40}$ and $q^{30}$, respectively, \eqref{Lovejoy-Osburn-2010} takes the form
\begin{align*}
&\sum_{n=-\infty}^\infty\frac{(-1)^nq^{25n^2+20n}}{1+q^{50n+10}}
-\sum_{n=-\infty}^\infty\frac{(-1)^nq^{25n^2+80n+60}}{1+q^{50n+70}}
\\
&\quad=-\frac{2q^{25}(q^{10},q^{40};q^{50})_\infty(-q^{25};q^{25})_\infty^2}
{(q^{20},q^{30};q^{50})_\infty(-q^{5},-q^{20};q^{25})_\infty}
\sum_{n=-\infty}^\infty\frac{(-1)^nq^{25n^2+50n}}{1+q^{50n+40}}
\\
&\ \ \
\quad\quad+\frac{(q^{10},q^{15};q^{25})_\infty(q^{10};q^{10})_\infty(q^{100};q^{100})_\infty}
{(-q^5,-q^{20};q^{25})_\infty(-q^{10};q^{10})_\infty(-q^{10},-q^{40};q^{50})_\infty}.
\end{align*}
Replacing $n$ by $n-1$ in second summation on the left side of the above relation, we find that
\begin{align*}
&\sum_{n=-\infty}^\infty\frac{(-1)^nq^{25n^2+20n}}{1+q^{50n+10}}
-\sum_{n=-\infty}^\infty\frac{(-1)^nq^{25n^2+30n+5}}{1+q^{50n+20}}
\\
&\quad=-\frac{2qJ_{5,50}}{J_{25,50}}\sum_{n=-\infty}^\infty\frac{(-1)^nq^{25n^2+50n+24}}{1+q^{50n+40}}
+U_2,
\end{align*}
that is,
\begin{align}\label{main-result-10-claim-2-M2}
&P_2-P_1=\frac{2q^4J_{5,50}}{J_{25,50}}P_4-q^3 U_2.
\end{align}
Plugging \eqref{main-result-10-claim-1-M2} and \eqref{main-result-10-claim-2-M2} into \eqref{main-result-10-disse-1-M2}, we obtain that
\begin{align}\label{10-eq---2}
\sum_{n=-\infty}^\infty\frac{(-1)^nq^{n^2+2n}}{1+q^{10n}}
=\left(1-\frac{2qJ_{15,50}}{J_{25,50}}+\frac{2q^4J_{5,50}}{J_{25,50}}\right)P_4
+\frac{1}{2}U_1-q^3 U_2.
\end{align}
Employing the following  5-dissection formula \cite[Lemma 3.1]{Lovejoy-Osburn-2008}
\begin{align}\label{Lovejoy-Osburn-2008-mod10}
\frac{(q;q)_\infty}{(-q;q)_\infty}&=\frac{(q^{25};q^{25})_\infty}{(-q^{25};q^{25})_\infty}
-2q(q^{15},q^{35},q^{50};q^{50})_\infty+2q^4(q^5,q^{45},q^{50};q^{50})_\infty
\nonumber\\
&=J_{25,50}-2qJ_{15,50}+2q^4J_{5,50}
\end{align}
and using \eqref{10-eq---2}, we infer that
\begin{align*}
\sum_{n=-\infty}^\infty\frac{(-1)^nq^{n^2+2n}}{1+q^{10n}}
=\frac{(q;q)_\infty}{(-q;q)_\infty}\cdot\frac{P_4}{J_{25,50}}+\frac{1}{2}U_1-q^3 U_2.
\end{align*}
This completes the proof.\qed

\begin{lem}\label{mod10-eq2-lem-M2}
Let
\begin{align*}
&V_1=\frac{J_{5,100}^2J_{10,100}J_{20,100}^3J_{30,100}^4J_{45,100}^2J_{50,100}^2}{J_{100}^{12}},
\\
&V_2=\frac{J_{5,100}J_{10,100}^3J_{15,100}J_{20,100}^2J_{30,100}^3J_{35,100}J_{40,100}
J_{45,100}J_{50,100}}{J_{100}^{12}}.
\end{align*}
We have
\begin{align}\label{mod10-M2-eq2}
\sum_{n=-\infty}^\infty
\frac{(-1)^nq^{n^2+4n}}{1+q^{10n}}=\frac{1}{2}V_1-q^7V_2-
\frac{(q;q)_\infty}{(-q;q)_\infty}\cdot\frac{q^{21}}{J_{25,50}}
\sum_{n=-\infty}^\infty\frac{(-1)^nq^{25n^2+50n}}{1+q^{50n+30}}.
\end{align}
\end{lem}

\pf
To prove \eqref{mod10-M2-eq2}, we rewrite the left side of \eqref{mod10-M2-eq2} as follows:
\begin{align}\label{main-result-10-disse-2-M2}
\sum_{n=-\infty}^\infty\frac{(-1)^nq^{n^2+4n}}{1+q^{10n}}
&=\sum_{n=-\infty}^\infty\frac{(-1)^nq^{25n^2+20n}}{1+q^{50n}}
-\sum_{n=-\infty}^\infty\frac{(-1)^nq^{25n^2+30n+5}}{1+q^{50n+10}}
+\sum_{n=-\infty}^\infty\frac{(-1)^nq^{25n^2+40n+12}}{1+q^{50n+20}}
\nonumber\\
&\quad\quad-\sum_{n=-\infty}^\infty\frac{(-1)^nq^{25n^2+50n+21}}{1+q^{50n+30}}
+\sum_{n=-\infty}^\infty\frac{(-1)^nq^{25n^2+60n+32}}{1+q^{50n+40}}.
\end{align}
Let $S_0, S_1, S_2, S_3$ and $S_4$ denote the summations on the right hand side
of \eqref{main-result-10-disse-2-M2}, respectively.
We claim that the following two relations hold:
\begin{align}\label{main-result-10-claim-3-M2}
&S_0-S_1=-\frac{2q^4J_{5,50}}{J_{25,50}}S_3+\frac{1}{2}V_1,
\\[3pt]
\label{main-result-10-claim-4-M2}
&S_2+S_4=\frac{2qJ_{15,50}}{J_{25,50}}S_3-q^7V_2.
\end{align}
 These two relations turn out to be specializations of the following  identity
  due to Chan \cite[Theorem 2.1]{Chan}:
\begin{align}\label{Chan}
&\frac{(a_1,q/a_1,\ldots, a_r,q/a_r;q)_\infty(q;q)_\infty^2}{(b_1,q/b_1,\ldots, b_s,q/b_s;q)_\infty}
\nonumber\\&=\frac{(a_1/b_1,qb_1/a_1,\ldots, a_r/b_1,qb_1/a_r;q)_\infty}{(b_2/b_1,qb_1/b_2,\ldots ,b_s/b_1,qb_1/b_s;q)_\infty}
\sum_{n=-\infty}^\infty\frac{(-1)^{(s-r)n}q^{(s-r)n(n+1)/2}}{1-b_1q^n}\left(
\frac{a_1\cdots a_rb_1^{s-r-1}}{b_2\cdots b_s}\right)^n
\nonumber\\
&\quad+\mathrm{idem}(b_1;b_2,\ldots,b_s),
\end{align}
where $\mathrm{idem}(b_1;b_2,\ldots,b_s)$ is defined by the relation
\begin{align}\label{Chan-idem}
&F(b_1,b_2,\cdots,b_m)+\mathrm{idem}(b_1;b_2,\cdots,b_m)
\nonumber\\[3pt]
&\quad=F(b_1,b_2,\cdots,b_m)+F(b_2,b_1,b_3,\cdots,b_m)
+\cdots+F(b_m,b_2,\cdots,b_{m-1},b_1).
\end{align}

To derive  \eqref{main-result-10-claim-3-M2} from \eqref{Chan}, we set $r=2,s=3$ and replace $q,~a_1,~a_2,~b_1,~b_2,~b_3$  by $q^{50}, \ q^{30},\ q^{5},\ -1,\ -q^{10},\ -q^{30}$, respectively, to get
\begin{align*}
&\sum_{n=-\infty}^\infty(-1)^nq^{25n^2}
\left(\frac{q^{20n}}{1+q^{50n}}-\frac{q^{30n+5}}{1+q^{50n+10}}\right)
\\&\quad=-\frac{q^4(-1,q^{10},-q^{25},-q^{25},q^{40},-q^{50};q^{50})_\infty}
{(-q^{5},q^{20},-q^{20},q^{30},-q^{30},-q^{45};q^{50})_\infty}
\sum_{n=-\infty}^\infty\frac{(-1)^nq^{25n^2+50n+21}}{1+q^{25n+30}}
\\&
\quad\quad\quad  +\frac{(q^5,q^{10},q^{40},q^{45};q^{50})_\infty(q^{20},q^{30},q^{50};q^{50})_\infty^2}
{(-1,-q^{5},-q^{10},-q^{40},-q^{45},-q^{50};q^{50})_\infty(-q^{20},-q^{30};q^{50})_\infty^2}\\
&\quad=-\frac{2q^4J_{5,50}}{J_{25,50}}\sum_{n=-\infty}^\infty\frac{(-1)^nq^{25n^2+50n+21}}{1+q^{25n+30}}
+V_1,
\end{align*}
which coincides with \eqref{main-result-10-claim-3-M2}.

The identity \eqref{main-result-10-claim-4-M2} can be justified by using the identity \eqref{Chan},
setting $r=2,s=3$ and replacing $q,~a_1,~a_2,~b_1,~b_2,~b_3$ by $q^{50}, \ q^{30},\ q^{55},\ -q^{20},\ -q^{40},\ -q^{30}$, respectively, we find that
\begin{align*}
&\sum_{n=-\infty}^\infty(-1)^nq^{25n^2}
\left(\frac{q^{40n}}{1+q^{50n+20}}+\frac{q^{60n+20}}{1+q^{50n+40}}\right)
\\&\quad=\frac{2(q^{20},q^{30};q^{50})_\infty(-q^{25};q^{25})_\infty^2}
{(-q^{10},-q^{15},-q^{35},-q^{40};q^{50})_\infty(q^{10},q^{40};q^{50})_\infty}
\sum_{n=-\infty}^\infty\frac{(-1)^nq^{25n^2+50n+10}}{1+q^{50n+30}}
\\&
\quad\quad\quad  -\frac{q^{-5}(q^5,q^{10},q^{40},q^{45};q^{50})_\infty(q^{20},q^{30},q^{50};q^{50})_\infty^2}
{(-q^{15},-q^{35};q^{50})_\infty(-q^{10},-q^{20},-q^{30},-q^{40};q^{50})_\infty^2}\\
&\quad=\frac{2qJ_{15,50}}{J_{25,50}}\sum_{n=-\infty}^\infty\frac{(-1)^nq^{25n^2+50n+9}}{1+q^{50n+30}}
-q^{-5}V_2,
\end{align*}
which gives \eqref{main-result-10-claim-4-M2}.

Now, we restate \eqref{mod10-M2-eq2} as follows:
\begin{align}\label{claim}
\sum_{n=-\infty}^\infty
\frac{(-1)^nq^{n^2+4n}}{1+q^{10n}}=-
\frac{(q;q)_\infty}{(-q;q)_\infty}\cdot\frac{q^{21}}{J_{25,50}}
S_3+\frac{1}{2}V_1-q^7V_2.
\end{align}
Substituting \eqref{main-result-10-claim-3-M2} and \eqref{main-result-10-claim-4-M2} into \eqref{main-result-10-disse-2-M2}, we deduce that
\begin{align*}
\sum_{n=-\infty}^\infty\frac{(-1)^nq^{n^2+4n}}{1+q^{10n}}
&=-\left(1-\frac{2qJ_{15,50}}{J_{25,50}}+\frac{2q^4J_{5,50}}{J_{25,50}}\right)S_3
+\frac{1}{2}V_1-q^7V_2.
\end{align*}
In view of the $5$-dissection formula \eqref{Lovejoy-Osburn-2008-mod10},
the above relation leads to \eqref{claim}.
This completes the proof.\qed

The following two lemmas are needed in the proofs of Theorems \ref{main-result-10-1-M2} and \ref{main-result-10-2-M2}.

\begin{lem}\label{U1-2q^3U2}
Assume that $U_1$ and $U_2$ are defined as in  Lemma \ref{mod10-lem-M2}.
Then the following identity holds:
\begin{align}\label{U1-2q^3U2-eq}
&\frac{1}{2}U_1-q^3U_2
\nonumber\\[5pt]
&\quad=\frac{(q;q)_\infty}{(-q;q)_\infty}\cdot\frac{J_{50}^4J_{25,50}}
{J_{5,10}^6}
\left\{A_0(q^5)+qA_1(q^5)+q^2 A_2(q^5)+q^3 A_3(q^5)+q^4 A_4(q^5)\right\},
\end{align}
where $A_0(q),\,A_1(q),\,A_2(q),\,A_3(q),\,A_4(q)$  are defined as in Theorem \ref{main-result-10-1-M2}.
\end{lem}

\begin{lem}\label{V1-U1+2q^3U2-2q^7V2}
Recall that $U_1$ and $U_2$ are defined in Lemma \ref{mod10-lem-M2} and $V_1$ and $V_2$ are defined as in Lemma \ref{mod10-eq2-lem-M2}. The following identity holds:
\begin{align}\label{V1-U1+2q^3U2-2q^7V2-eq}
&\frac{1}{2}V_1-\frac{1}{2}U_1+q^3U_2-q^7V_2
\nonumber\\[5pt]
&\quad=\frac{(q;q)_\infty}{(-q;q)_\infty}\cdot\frac{J_{50}^4J_{25,50}}{J_{5,10}^6}
\left\{B_0(q^5)+qB_1(q^5)+q^2 B_2(q^5)+q^3 B_3(q^5)+q^4 B_4(q^5)\right\},
\end{align}
where $B_0(q),\,B_1(q),\,B_2(q),\,B_3(q),\,B_4(q)$ are defined as in Theorem \ref{main-result-10-2-M2}.
\end{lem}

To present the proofs of Lemmas \ref{U1-2q^3U2} and \ref{V1-U1+2q^3U2-2q^7V2}, we
shall use techniques in the theory of modular forms. Recall that the Dedekind $\eta$-function is defined by
\[\eta(\tau)=q^{\frac{1}{24}}(q;q)_\infty,\]
where $\tau\in \mathcal{H}=\{\tau\in\mathcal{C}:\mathrm{Im} \tau>0\}$ and $q=e^{2\pi i \tau}$, and the generalized Dedekind $\eta$-function is defined by
\begin{align}\label{eta-defi}
\eta_{\delta,g}(\tau)=q^{P(g/\delta)\delta/2}\prod_{\mbox{\tiny$\begin{array}{c}
n>0\\
n\equiv g\pmod\delta\end{array}$}} (1-q^n)
\prod_{\mbox{\tiny$\begin{array}{c}
n>0\\
n\equiv -g\pmod\delta\end{array}$}} (1-q^n),
\end{align}
 where $g,\delta \in \mathbb{Z}^+$,  $0<g<\delta$,  $P(t)=\{t\}^2-\{t\}+\frac{1}{6}$ is the second Bernoulli function, and $\{t\}=t-[t]$ is the fractional parts of $t$. In
  particular,  \[\eta_{\delta,0}(\tau)=q^{\frac{\delta}{12}}(q^{\delta};q^{\delta})_\infty^2\]
  and
\[\eta_{\delta,\frac{\delta}{2}}(\tau)
=q^{-\frac{\delta}{24}}(q^{\frac{\delta}{2}};q^{\delta})_\infty^2.\]

Let $N$ be a fixed positive integer. A generalized  Dedekind $\eta$-quotient of level $N$ has the form
\begin{align}\label{eta-quotient}
f(\tau)=\prod_{\delta |N \atop 0<g<\delta} \eta_{\delta,g}^{r_{\delta,g}}(\tau),
\end{align}
where
\[r_{\delta,g}\in \left\{\begin{array}{ll}
\frac{1}{2}\mathbb{Z}, &\text{if}~g=0~\text{or}~g=\delta/2,\\[5pt]
\mathbb{Z}, &\text{otherwise.}
\end{array}
\right.
\]

 Suppose that $f$ is a modular function with respect to the congruence subgroup $\Gamma$ of $\Gamma_0(1)$.
 For $A\in\Gamma_0(1)$, there exits a cusp given by $\zeta=A^{-1}\infty$.
If
\[f(A^{-1}\tau)=\sum_{m=m_0}^\infty b_mq^{m/N}\]
and $b_{m_0}\neq0$, then we say that $m_0$ is the order of $f$ at $\zeta$ with respect to $\Gamma$, and it is denoted by $ORD(f,\zeta,\Gamma)$.

Robins \cite{Robins-1994} found sufficient conditions for a generalized $\eta$-quotient
to be a modular function on $\Gamma_1(N)$.
\begin{thm}\label{Robins-1994}
Let $f(\tau)$ be a generalized $\eta$-quotient defined as in \eqref{eta-quotient}. If
\begin{align*}
\sum_{\delta |N\atop 0\leq g\leq\delta}\delta P_2\left(\frac{g}{\delta}\right)
r_{\delta,g}\equiv0\pmod2,
\end{align*}
and
\begin{align*}
\sum_{\delta |N\atop 0\leq g\leq\delta}\frac{N}{\delta} P_2(0)
r_{\delta,g}\equiv0\pmod2,
\end{align*}
then $f(\tau)$ is a modular function on $\Gamma_1(N)$.
\end{thm}

The following theorem of Garvan and Liang \cite{Garvan-Liang-2016} can be used to prove generalized $\eta$-quotient identities.

\begin{thm}\label{eta-lem}
Let $f_1(\tau),f_2(\tau),\ldots,f_n(\tau)$ be generalized $\eta$-quotients that are modular functions on $\Gamma_1(N)$. Let $\mathcal{S}_N$ be a set of inequivalent cusps for $\Gamma_1(N)$. Define the constant
\begin{eqnarray}\label{constantB}
B=\sum_{\substack{s\in\mathcal{S}_N \\ s\neq i\infty}}
\text{min}(\{ORD(f_j,s,\Gamma_1(N)):1\leq j\leq n\}\cup\{0\}).
\end{eqnarray}
Let
\begin{eqnarray}
g(\tau)=\alpha_1f_1(\tau)+\alpha_2f_2(\tau)+\cdots+\alpha_nf_n(\tau)+1,
\end{eqnarray}
where each $\alpha_j\in\mathbb{C}$. Then
\[g(\tau)\equiv0\]
if and only if
\begin{eqnarray}
ORD(g(\tau),i\infty,\Gamma_1(N))>-B.
\end{eqnarray}
\end{thm}

We are now in a position to  prove Lemmas  \ref{U1-2q^3U2} and \ref{V1-U1+2q^3U2-2q^7V2} with the aid of  Theorem \ref{eta-lem}.

\vskip 0.2cm
\noindent{\it Proof of Lemma \ref{U1-2q^3U2}.}
The relation \eqref{U1-2q^3U2-eq} can be stated as follows by the definitions of $A_0(q)$, $A_1(q)$, $A_2(q)$, $A_3(q)$, $A_4(q)$ as given in
\eqref{10-eq-A0}, \eqref{10-eq-A1}, \eqref{10-eq-A2}, \eqref{10-eq-A3} and \eqref{10-eq-A4}:
\begin{footnotesize}
\begin{align}\label{generfun-1}
&\frac{1}{2}\frac{J_{10,100}^4J_{15,100}^2J_{30,100}J_{35,100}^2J_{40,100}^3J_{50,100}^2}
{J_{100}^{12}}
-q^3\frac{J_{5,100}J_{10,100}^3J_{15,100}J_{20,100}J_{30,100}^3
J_{35,100}J_{40,100}^2J_{45,100}J_{50,100}}{J_{100}^{12}}
\nonumber\\[5pt]
&=\frac{(q;q)_\infty}{(-q;q)_\infty}\cdot\frac{J_{50}^4J_{25,50}}{J_{5,10}^6}
\left\{\frac{J_{10}^4J_{15,100}^2J_{25,100}^8J_{35,100}^2}{2J_{20,100}^4J_{30,100}^3
J_{40,100}J_{100}^6}
+\frac{2q^5 J_{5}^3J_{50}J_{25,100}J_{100}^5}{J_{5,100}^2J_{10,100}J_{40,100}^3
J_{45,100}^2}\right.
\nonumber\\[5pt]&\ \ \ \ \ \ \ \
+\frac{8q^{10}J_{5}^2J_{10}^4 J_{10,100}J_{15,100}J_{35,100}}{J_{20,100}^2
J_{30,100}J_{40,100}J_{50,100}^2J_{100}}
-\frac{16q^{10}J_{5}^2 J_{10,100}J_{15,100}^2J_{30,100}J_{35,100}^2}{J_{20,100}
J_{25,100}^2J_{100}^3}
\nonumber\\[5pt]&
+q\left(\frac{J_{10}^4 J_{15,100}^3J_{25,100}^6J_{35,100}^3}
{J_{20,100}^4J_{30,100}^3J_{40,100}J_{100}^6}
+\frac{8q^5J_{5}^2J_{50}^2J_{10,100}^2J_{15,100}^2J_{35,100}^2J_{40,100}}
{J_{5,100}J_{20,100}^2J_{30,100}J_{45,100}J_{100}^4}
-\frac{8q^5J_{5}^2 J_{10,100}^4J_{30,100}^2J_{40,100}^3}
{J_{5,100}^2J_{15,100}J_{35,100}J_{45,100}^2J_{100}^3}\right.
\nonumber\\[5pt]&\ \ \ \ \ \ \ \
+\frac{4q^{10}J_{5}J_{10}^2 J_{5,100}^2J_{25,100}^3J_{45,100}^2}
{J_{20,100}^2J_{40,100}J_{100}^4J_{50}}
\left.-\frac{8q^{15}J_{5}^2J_{5,100}^2J_{20,100}J_{30,100}J_{45,100}^2J_{50,100}}
{J_{25,100}^2J_{40,100}^2J_{100}^3}\right)
\nonumber\\[5pt]&
+q^2\left(\frac{2J_{10}^4J_{15,100}^4J_{25,100}^4J_{35,100}^4}
{J_{20,100}^4J_{30,100}^3J_{40,100}J_{100}^6}
+\frac{2q^5J_{5}^2J_{10,100}J_{25,100}^4J_{30,100}}
{J_{15,100}J_{20,100}J_{35,100}J_{100}^3}
-\frac{8q^5J_{5}^2J_{50}^2 J_{15,100}^3J_{20,100}J_{30,100}J_{35,100}^3}
{J_{5,100}J_{25,100}^2J_{40,100}^2J_{45,100}J_{100}^4}\right.
\nonumber\\[5pt]&\ \ \ \ \ \ \ \
+\frac{4q^{10}J_{5}^2J_{50}^2 J_{5,100}J_{20,100}J_{30,100}J_{45,100}}
{J_{40,100}^2J_{100}^4}
\left.+\frac{8q^{10}J_{5}^2 J_{5,100}J_{10,100}J_{30,100}J_{45,100}}
{J_{20,100}J_{100}^3}\right)
\nonumber\\[5pt]&
+q^3\left(-\frac{J_{5}J_{10}^2J_{25,100}^7}
{J_{50}J_{20,100}^2J_{40,100}J_{100}^4}
+\frac{4J_{10}^4J_{15,100}^5J_{25,100}^2J_{35,100}^5}
{J_{20,100}^4J_{30,100}^3J_{40,100}J_{100}^6}
+\frac{2q^5 J_{5}^2J_{50}^2 J_{10,100}^2J_{25,100}^2J_{40,100}}
{J_{20,100}^2J_{30,100}J_{100}^4}\right.
\nonumber\\[5pt]&\ \ \ \ \ \ \ \
-\frac{12q^5 J_{5}^2 J_{10,100}J_{25,100}^2J_{30,100}}
{J_{20,100}J_{100}^3}
\left.-\frac{16q^{10}J_{5}^3J_{100}^6}{J_{20,100}^2J_{25,100}^3J_{40,100}J_{50}}
+\frac{8q^{10}J_{5}^3J_{50}J_{10,100}J_{100}^5}
{J_{20,100}^3J_{25,100}^3J_{30,100}^2}\right)
\nonumber\\[5pt]&
+q^4\left(-\frac{2J_{5}^2J_{10,100}J_{25,100}^4J_{30,100}}
{J_{5,100}J_{20,100}J_{45,100}J_{100}^3}
-\frac{J_{5}^2J_{50}^2 J_{10,100}^2J_{25,100}^4J_{40,100}}
{J_{5,100}J_{20,100}^2J_{30,100}J_{45,100}J_{100}^4}
+\frac{8J_{10}^4J_{15,100}^6J_{35,100}^6}{J_{20,100}^4J_{30,100}^3J_{40,100}J_{100}^6}\right.
\nonumber\\[5pt]&\ \ \ \ \ \ \ \
-\frac{16q^5J_{5}^2J_{10,100}J_{15,100}J_{30,100}J_{35,100}}
{J_{20,100}J_{100}^3}
-\frac{8q^5J_{5}^3J_{50}J_{100}^5J_{10,100}}{J_{5,100}J_{20,100}^3J_{25,100}
J_{30,100}^2J_{45,100}}
\nonumber\\[5pt]&\ \ \ \ \ \ \ \
\left.\left.-\frac{16q^{15}J_{5}J_{10}^2J_{5,100}^4J_{45,100}^4}
{J_{50}J_{100}^4J_{20,100}^2J_{25,100}J_{40,100}}\right)\right\}.
\end{align}
\end{footnotesize}
Multiplying both sides by $2J_{10,100}^{-4}J_{15,100}^{-2}J_{30,100}^{-1}J_{35,100}^{-2}
J_{40,100}^{-3}J_{50,100}^{-2}J_{100}^{12}$, the above relation becomes
\begin{footnotesize}
\begin{align*}
&1-\frac{2\eta_{100,5}(\tau)\eta_{100,20}(\tau)
\eta_{100,30}(\tau)^2\eta_{100,45}(\tau)}
{\eta_{100,10}(\tau)\eta_{100,15}(\tau)
\eta_{100,35}(\tau)\eta_{100,40}(\tau)\eta_{100,50}(\tau)}
=\frac{\eta_{1,0}(\tau)\eta_{50,0}(\tau)^{\frac{5}{2}}\eta_{50,25}(\tau)}
{\eta_{2,0}(\tau)^{\frac{1}{2}}\eta_{10,0}(\tau)^3\eta_{10, 5}(\tau)^6}
\Bigg(\eta_{100,25}(\tau)^8
\\[5pt]
+&\frac{4\eta_{10,5}(\tau)^{\frac{1}{2}}\eta_{100,20}(\tau)^3\eta_{100,25}(\tau)^3\eta_{100,30}(\tau)^2}
{\eta_{100,10}(\tau)^2\eta_{100,40}(\tau)^3}
+\frac{16\eta_{10,0}(\tau)\eta_{10,5}(\tau)\eta_{100,10}(\tau)\eta_{100,20}(\tau)^2\eta_{100,30}(\tau)^2}
{\eta_{100,0}(\tau)\eta_{100,15}(\tau)\eta_{100,35}(\tau)\eta_{100,50}(\tau)^2}
\\[5pt]
-&\frac{32\eta_{10,5}(\tau)\eta_{100,20}(\tau)\eta_{100,30}(\tau)^2}
{\eta_{100,10}(\tau)\eta_{100,25}(\tau)^2\eta_{100,40}(\tau)\eta_{100,50}(\tau)}
+2\eta_{100,15}(\tau)\eta_{100,25}(\tau)^6\eta_{100,35}(\tau)
+\frac{16\eta_{10,5}(\tau)}{\eta_{100,5}(\tau)\eta_{100,45}(\tau)}
\\[5pt]
-&\frac{16\eta_{10,0}(\tau)\eta_{100,25}(\tau)^2\eta_{100,30}(\tau)}
{\eta_{100,0}(\tau)\eta_{100,15}(\tau)\eta_{100,35}(\tau)\eta_{100,50}(\tau)^2}
+\frac{8\eta_{100,5}(\tau)^3\eta_{100,20}(\tau)\eta_{100,25}(\tau)^4
\eta_{100,30}(\tau)^2\eta_{100,45}(\tau)^3}
{\eta_{100,10}(\tau)\eta_{100,15}(\tau)\eta_{100,35}(\tau)\eta_{100,40}(\tau)\eta_{100, 50}(\tau)}
\\[5pt]
-&\frac{16\eta_{100,5}(\tau)^4\eta_{100,20}(\tau)^3\eta_{100,30}(\tau)^2
\eta_{100,45}(\tau)^4}{\eta_{100,10}(\tau)^2\eta_{100,40}(\tau)^3}
+4\eta_{100,15}(\tau)^2\eta_{100,25}(\tau)^4\eta_{100,35}(\tau)^2
\\[5pt]
+&\frac{4\eta_{100,5}(\tau)^2\eta_{100,20}(\tau)\eta_{100,25}(\tau)^6
\eta_{100,30}(\tau)^2\eta_{100,45}(\tau)^2}
{\eta_{100,10}(\tau)\eta_{100,15}(\tau)\eta_{100,35}(\tau)
\eta_{100, 40}(\tau)\eta_{100, 50}(\tau)}
-\frac{16\eta_{10,5}(\tau)^{\frac{1}{2}}\eta_{100,15}(\tau)^2\eta_{100,20}(\tau)^3\eta_{100,30}(\tau)^2\eta_{100,35}(\tau)^2}
{\eta_{100,10}(\tau)^2\eta_{100,25}(\tau)\eta_{100,40}(\tau)^3}
\\[5pt]
+&\frac{8\eta_{100,5}(\tau)^3\eta_{100,20}(\tau)^3\eta_{100,25}(\tau)^2
\eta_{100,30}(\tau)^2\eta_{100,45}(\tau)^3}
{\eta_{100,10}(\tau)^2\eta_{100,40}(\tau)^3}
+\frac{16\eta_{100,5}(\tau)^3\eta_{100,20}(\tau)\eta_{100,25}(\tau)^2
\eta_{100,30}(\tau)^2\eta_{100,45}(\tau)^3}
{\eta_{100,10}(\tau)\eta_{100,40}(\tau)\eta_{100,50}(\tau)}
\\[5pt]
-&\frac{2\eta_{100,5}(\tau)\eta_{100,20}(\tau)\eta_{100,25}(\tau)^8
\eta_{100,30}(\tau)^2\eta_{100,45}(\tau)}{\eta_{100,10}(\tau)\eta_{100,15}(\tau)
\eta_{100,35}(\tau)\eta_{100,40}(\tau)\eta_{100,50}(\tau)}
+8\eta_{100,15}(\tau)^3\eta_{100,25}(\tau)^2\eta_{100,35}(\tau)^3
\\[5pt]
+&4\eta_{100,5}(\tau)^2\eta_{100,25}(\tau)^4\eta_{100,45}(\tau)^2
-\frac{24\eta_{100,5}(\tau)^2\eta_{100,20}(\tau)\eta_{100,25}(\tau)^4
\eta_{100,30}(\tau)^2\eta_{100,45}(\tau)^2}
{\eta_{100,10}(\tau)\eta_{100,40}(\tau)\eta_{100,50}(\tau)}
\\[5pt]
-&\frac{32\eta_{10,5}(\tau)^{\frac{1}{2}}\eta_{100,5}(\tau)^2\eta_{100,20}(\tau)\eta_{100,30}(\tau)^2\eta_{100,45}(\tau)^2}
{\eta_{100,10}(\tau)\eta_{100,25}(\tau)\eta_{100,40}(\tau)\eta_{100,50}(\tau)}
\\[5pt]
+&\frac{16\eta_{10,5}(\tau)^{\frac{1}{2}}\eta_{100,5}(\tau)^2\eta_{100,45}(\tau)^2}{\eta_{100,25}(\tau)}
-\frac{4\eta_{100,5}(\tau)\eta_{100,20}(\tau)\eta_{100,25}(\tau)^6
\eta_{100,30}(\tau)^2\eta_{100,45}(\tau)}
{\eta_{100,10}(\tau)\eta_{100,40}(\tau)\eta_{100,50}(\tau)}
\\[5pt]
-&2\eta_{100, 5}(\tau)\eta_{100, 25}(\tau)^6\eta_{100, 45}(\tau)
+16\eta_{100, 15}(\tau)^4\eta_{100, 35}(\tau)^4
\\[5pt]
-&\frac{32\eta_{10,5}(\tau)\eta_{100,20}(\tau)\eta_{100,30}(\tau)^2}
{\eta_{100,10}(\tau)\eta_{100,15}(\tau)\eta_{100,35}(\tau)\eta_{100,40}(\tau)\eta_{100,50}(\tau)}
-\frac{16\eta_{10,5}(\tau)}{\eta_{100,15}(\tau)\eta_{100,35}(\tau)}
\\[5pt]
-&\frac{32\eta_{100, 5}(\tau)^5\eta_{100, 20}(\tau)\eta_{100, 30}(\tau)^2
\eta_{100, 45}(\tau)^5}{\eta_{100, 10}(\tau)\eta_{100, 15}(\tau)\eta_{100, 35}(\tau)\eta_{100, 40}(\tau)\eta_{100, 50}(\tau)}\Bigg).
\end{align*}
\end{footnotesize}
By means of Theorem \ref{Robins-1994}, it can be shown that each term of the above identity is a modular function with respect to $\Gamma_1(100)$. Applying the algorithm  of Garvan and Liang \cite{Garvan-Liang-2016}, we find that the constant $B$  in \eqref{constantB} is equal to  $-600$. Thus, by Theorem \ref{eta-lem},
 \eqref{generfun-1} holds as long as one can verify it holds for
 the coefficients  of $q^{n}$ for $1\leq n \leq 600$. This task can be
 easily accomplished, and hence the proof is complete.\qed

\vskip 0.2cm
\noindent{\it Proof of Lemma  \ref{V1-U1+2q^3U2-2q^7V2}.}
By the definitions of $B_0(q)$, $B_1(q)$, $B_2(q)$, $B_3(q)$, $B_4(q)$ as given in \eqref{10-eq-B0}, \eqref{10-eq-B1}, \eqref{10-eq-B2}, \eqref{10-eq-B3} and \eqref{10-eq-B4}, we see that  \eqref{V1-U1+2q^3U2-2q^7V2-eq} can be expressed as
\begin{footnotesize}
\begin{align}\label{generfun-2}
&-\frac{1}{2}\frac{J_{10,100}^4J_{15,100}^2J_{30,100}J_{35,100}^2J_{40,100}^3J_{50,100}^2}
{J_{100}^{12}}
+q^3\frac{J_{5,100}J_{10,100}^3J_{15,100}J_{20,100}J_{30,100}^3
J_{35,100}J_{40,100}^2J_{45,100}J_{50,100}}{J_{100}^{12}}
\nonumber\\
&+\frac{1}{2}\frac{J_{5,100}^2J_{10,100}J_{20,100}^3J_{30,100}^4J_{45,100}^2J_{50,100}^2}{J_{100}^{12}}
-q^7\frac{J_{5,100}J_{10,100}^3J_{15,100}J_{20,100}^2J_{30,100}^3J_{35,100}J_{40,100}
J_{45,100}J_{50,100}}{J_{100}^{12}}
\nonumber\\
&=\frac{(q;q)_\infty}{(-q;q)_\infty}\cdot\frac{J_{50}^4J_{25,50}}{J_{5,10}^6}
\left\{-\frac{q^5J_{5}^5J_{5,100}J_{25,100}J_{45,100}}{J_{50}J_{100}J_{20,100}^2J_{30,100}J_{40,100}}
+\frac{4q^5J_{5}^4J_{100}^2J_{15,100}J_{35,100}}{J_{5,100}J_{20,100}^2J_{40,100}J_{45,100}J_{50,100}}\right.
\nonumber\\[5pt]&
\quad-\frac{16q^{10}J_{10}^{\frac{7}{2}}J_{5,10}^{\frac{1}{2}}J_{5,100}J_{10,100}J_{45,100}}{J_{100}J_{20,100}J_{25,100}^3}
+\frac{16q^{10}J_{5}^2J_{10,100}J_{15,100}^2J_{30,100}J_{35,100}^2}{J_{100}^3 J_{20,100}J_{25,100}^2}
-\frac{8q^{15}J_{5}J_{10}^2J_{5,100}^2J_{25,100}^3J_{45,100}^2}{J_{50}J_{100}^4J_{20,100}J_{40,100}^2}
\nonumber\\[5pt]&
\quad-\frac{4q^{15}J_{5}^4J_{100}^2J_{5,100}J_{45,100}}{J_{15,100}J_{20,100}J_{35,100}J_{40,100}^2J_{50,100}}
-\frac{16q^{20}J_{5,100}^2J_{10,100}J_{30,100}J_{45,100}^2J_{5}^2}{J_{100}^3J_{25,100}^2J_{40,100}}
\nonumber\\[5pt]&
+q\left(-\frac{2q^5J_{5}^2J_{50}^2J_{10,100}J_{25,100}^2}{J_{100}^4J_{20,100}}
+\frac{8q^5J_{5}^3J_{10,100}J_{15,100}^3J_{30,100}J_{35,100}^3}{J_{50}J_{20,100}J_{25,100}J_{100}^6}\right.
+\frac{4q^{10}J_{10}^4J_{5,100}^2J_{10,100}J_{30,100}J_{45,100}^2}{J_{100}^6J_{20,100}J_{50,100}}
\nonumber\\[5pt]&
\quad+\frac{2q^{10}J_{10}^4J_{5,100}J_{10,100}J_{25,100}^2J_{30,100}J_{45,100}}{J_{100}^6J_{40,100}J_{50,100}}
-\frac{16q^{10}J_{5}^3J_{50}J_{100}^5}{J_{20,100}^2J_{25,100}^3J_{30,100}J_{40,100}}
+\frac{12q^{15}J_{5}^2J_{5,100}^2J_{10,100}J_{30,100}J_{45,100}^2}{J_{100}^3J_{20,100}J_{25,100}^2}
\nonumber\\[5pt]&
\left.
\quad+\frac{14q^{15}J_{5}^2J_{5,100}J_{10,100}J_{30,100}J_{45,100}}{J_{100}^3J_{40,100}}
-\frac{16q^{20}J_{5}^2J_{50}^2J_{5,100}^4J_{10,100}J_{45,100}^4}{J_{100}^4J_{15,100}J_{20,100}J_{35,100}J_{25,100}^4}\right)
\nonumber\\&
+q^2\left(\frac{q^5J_{5}^3J_{10,100}J_{25,100}^5J_{30,100}}{J_{50}J_{100}^6J_{40,100}}
+\frac{4q^5J_{5}^2J_{50}^2J_{5,100}^2J_{20,100}^2J_{30,100}^3J_{45,100}^2}
{J_{100}^4J_{10,100}^2J_{15,100}J_{35,100}J_{40,100}^3}\right.
\nonumber\\[5pt]&
\quad -\frac{8q^5J_{10}^4J_{5,100}J_{15,100}^5J_{35,100}^5J_{45,100}}{J_{100}^6J_{20,100}^4J_{30,100}^3J_{40,100}}
+\frac{12q^5J_{5}J_{10}^2J_{15,100}^4J_{35,100}^4}{J_{50}J_{100}^4J_{25,100}J_{20,100}^2J_{40,100}}
+\frac{4q^{10}J_{5}J_{50}J_{10,100}^6J_{30,100}^2J_{40,100}^5}{J_{25,100}J_{100}^{12}}
\nonumber\\[5pt]&
\quad
-\frac{12q^{10}J_{5}^2J_{10,100}J_{25,100}^2J_{30,100}}{J_{100}^3J_{40,100}}
\nonumber\left.
-\frac{16q^{10}J_{5}^2J_{5,100}J_{10,100}J_{30,100}J_{45,100}}{J_{100}^3J_{20,100}}
-\frac{8q^{15}J_{5}J_{10}^2J_{50}J_{10,100}^4J_{30,100}}{J_{100}^5J_{20,100}J_{25,100}}\right)
\nonumber\\&
+q^3\left(\frac{J_{5}^3J_{10,100}J_{25,100}^5J_{30,100}}{J_{50}J_{100}^6J_{20,100}}
-\frac{8q^5J_{5}^2J_{50}^2J_{10,100}^3J_{15,100}^2J_{35,100}^2J_{40,100}^2}
{J_{100}^4J_{5,100}J_{20,100}^3J_{30,100}^2J_{45,100}}\right.
+\frac{12q^5J_{5}^2J_{10,100}J_{25,100}^2J_{30,100}}{J_{100}^3J_{20,100}}
\nonumber\\[5pt]&
\quad
+\frac{4q^{10}J_{10}^4J_{5,100}^6J_{25,100}^2J_{30,100}J_{45,100}^6}{J_{100}^6J_{10,100}^4J_{15,100}J_{35,100}J_{40,100}^5}
\left.-\frac{16q^{15}J_{5}^3J_{100}^6J_{15,100}J_{35,100}}{J_{50}J_{5,100}J_{25,100}^3J_{40,100}^3J_{45,100}}
-\frac{16q^{20}J_{5}J_{10}^2J_{5,100}^4J_{45,100}^4}{J_{50}J_{100}^4J_{50,100}J_{25,100}J_{40,100}^2}\right)
\nonumber\\&
+q^4\left(\frac{2J_{5}^2J_{10,100}J_{25,100}^4J_{30,100}}{J_{100}^3J_{5,100}J_{20,100}J_{45,100}}
-\frac{2q^5J_{5}^2J_{50}^2J_{10,100}^5J_{25,100}^2J_{30,100}^2J_{40,100}}{J_{100}^4J_{5,100}^3J_{15,100}J_{35,100}J_{45,100}^3}
+\frac{16q^{10}J_{5}^2J_{50}^2J_{5,100}^2J_{10,100}J_{45,100}^2}{J_{100}^4J_{20,100}J_{25,100}^2}\right.
\nonumber\\[5pt]&
\quad
-\frac{16q^{10}J_{5}^2J_{10,100}J_{15,100}^2J_{30,100}J_{35,100}^2}{J_{100}^3J_{25,100}^2J_{40,100}}
+\frac{8q^{15}J_{5}^2J_{5,100}^2J_{10,100}J_{30,100}J_{45,100}^2}{J_{100}^3J_{15,100}J_{35,100}J_{40,100}}
+\frac{16q^{15}J_{5}J_{10}^2J_{5,100}^4J_{45,100}^4}{J_{50}J_{100}^4J_{20,100}^2J_{25,100}J_{40,100}}
\nonumber\\[5pt]&
\left.\left.\quad+\frac{32q^{20}J_{5}^2J_{10,100}^3J_{15,100}^2J_{35,100}^2}{J_{50}^2J_{100}^2J_{20,100}J_{25,100}^2}
\right)\right\}.
\end{align}
\end{footnotesize}
Multiplying both sides of above relation by $2J_{10,100}^{-4}J_{15,100}^{-2}J_{30,100}^{-1}J_{35,100}^{-2}
J_{40,100}^{-3}J_{50,100}^{-2}J_{100}^{12}$, the above relation takes the form
\begin{footnotesize}
\begin{align*}
&-1+\frac{2\eta_{100,5}(\tau)\eta_{100,20}(\tau)\eta_{100,30}(\tau)^2\eta_{100,45}(\tau)}
{\eta_{100,10}(\tau)\eta_{100,15}(\tau)\eta_{100,35}(\tau)\eta_{100,40}(\tau)\eta_{100,50}(\tau)}
\\
&+\frac{\eta_{100,5}(\tau)^2\eta_{100,20}(\tau)^3\eta_{100,30}(\tau)^3\eta_{100,45}(\tau)^2}
{\eta_{100,10}^3(\tau)\eta_{100,15}^2(\tau)\eta_{100,35}^2(\tau)\eta_{100,40}^3(\tau)}
-\frac{2\eta_{100,5}(\tau)\eta_{100,20}^2(\tau)\eta_{100,30}^2(\tau)\eta_{100,45}(\tau)}
{\eta_{100,10}(\tau)\eta_{100,15}(\tau)\eta_{100,35}(\tau)\eta_{100,40}^2(\tau)\eta_{100,50}(\tau)}
\\
&=\frac{\eta_{1,0}(\tau)\eta_{50,0}(\tau)^{\frac{5}{2}}\eta_{50,25}(\tau)}
{\eta_{2,0}(\tau)^{\frac{1}{2}}\eta_{10,0}(\tau)^3\eta_{10,5}(\tau)^6}
\Bigg(-\frac{2\eta_{10,5}(\tau)^3\eta_{100,10}(\tau)\eta_{100,20}(\tau)^3\eta_{100,30}(\tau)^3\eta_{100,40}(\tau)}
{\eta_{100,15}(\tau)^3\eta_{100,35}(\tau)^3}
\\[5pt]
&+\frac{8\eta_{10,5}(\tau)^{\frac{3}{2}}\eta_{100,20}(\tau)^2\eta_{100,25}(\tau)\eta_{100,30}(\tau)^3}{\eta_{100,50}(\tau)}
-\frac{32\eta_{10,5}(\tau)\eta_{100,10}(\tau)\eta_{100,20}(\tau)^3\eta_{100,30}(\tau)^3
\eta_{100,40}(\tau)}
{\eta_{100,15}(\tau)^3\eta_{100,25}(\tau)^4\eta_{100,35}(\tau)^3}
\\[5pt]
&+\frac{32\eta_{10,5}(\tau)\eta_{100,20}(\tau)\eta_{100,30}(\tau)^2}
{\eta_{100,10}(\tau)\eta_{100,25}(\tau)^2\eta_{100,40}(\tau)\eta_{100,50}(\tau)}
-\frac{8\eta_{10,5}(\tau)^2\eta_{100,5}(\tau)\eta_{100,20}(\tau)^3\eta_{100,30}(\tau)^3\eta_{100,45}(\tau)}
{\eta_{100,15}(\tau)^3\eta_{100,35}(\tau)^3\eta_{100,40}(\tau)\eta_{100,50}(\tau)}
\\[5pt]
&-\frac{16\eta_{100,5}(\tau)^3\eta_{100, 20}(\tau)^2\eta_{100, 25}(\tau)^4\eta_{100, 30}(\tau)^2\eta_{100, 45}(\tau)^3}
{\eta_{100, 10}(\tau)\eta_{100, 15}(\tau)\eta_{100, 35}(\tau)\eta_{100, 40}(\tau)^2\eta_{100, 50}(\tau)}
-\frac{32\eta_{10,0}(\tau)\eta_{100,5}(\tau)^4\eta_{100,45}(\tau)^4}
{\eta_{100,0}(\tau)\eta_{100,10}(\tau)^3\eta_{100,40}(\tau)^4\eta_{100,50}(\tau)^2}
\\[5pt]
&-\frac{4\eta_{10,5}(\tau)\eta_{100,20}(\tau)\eta_{100,25}(\tau)^2\eta_{100,30}(\tau)}
{\eta_{100,10}(\tau)\eta_{100,15}(\tau)^2\eta_{100,35}(\tau)^2\eta_{100,40}(\tau)}
+\frac{16\eta_{10,5}(\tau)^2\eta_{100,20}(\tau)^2\eta_{100,30}(\tau)^3}
{\eta_{100,5}(\tau)\eta_{100,25}(\tau)^2\eta_{100,45}(\tau)\eta_{100,50}(\tau)}
\\[5pt]
&+\frac{4\eta_{10,5}(\tau)^{\frac{1}{2}}\eta_{100,10}(\tau)\eta_{100,20}(\tau)^4\eta_{100,25}(\tau)\eta_{100,30}(\tau)^4}
{\eta_{100,15}(\tau)^3\eta_{100,35}(\tau)^3\eta_{100,50}(\tau)}
+\frac{8\eta_{10,5}(\tau)\eta_{100,10}(\tau)\eta_{100,20}(\tau)^3\eta_{100,30}(\tau)^4\eta_{100,40}(\tau)}
{\eta_{100,15}(\tau)^4\eta_{100,25}(\tau)^2\eta_{100,35}(\tau)^4\eta_{100,50}(\tau)}
\\[5pt]
&-\frac{32\eta_{10,5}(\tau)^{\frac{1}{2}}\eta_{100,5}(\tau)^2\eta_{100,20}(\tau)\eta_{100,30}(\tau)\eta_{100,45}(\tau)^2}
{\eta_{100,10}(\tau)\eta_{100,25}(\tau)\eta_{100,40}(\tau)}
+\frac{24\eta_{100, 5}(\tau)^4\eta_{100, 20}(\tau)\eta_{100, 30}(\tau)^2\eta_{100, 45}(\tau)^4}
{\eta_{100, 10}(\tau)\eta_{100, 40}(\tau)\eta_{100, 50}(\tau)}
\\[5pt]
&+\frac{28\eta_{100, 5}(\tau)^3\eta_{100, 20}(\tau)^2\eta_{100, 25}(\tau)^2\eta_{100, 30}(\tau)^2\eta_{100, 45}(\tau)^3}
{\eta_{100, 10}(\tau)\eta_{100, 40}(\tau)^2\eta_{100, 50}(\tau)}
-\frac{32\eta_{100, 5}(\tau)^6\eta_{100, 20}(\tau)\eta_{100, 30}(\tau)\eta_{100, 45}(\tau)^6}
{\eta_{100, 10}(\tau)\eta_{100, 25}(\tau)^2\eta_{100, 15}(\tau)\eta_{100, 35}(\tau)\eta_{100, 40}(\tau)}
\\[5pt]
&+\frac{2\eta_{10,5}(\tau)^{\frac{3}{2}}\eta_{100,20}(\tau)^3\eta_{100,25}(\tau)^5\eta_{100,30}(\tau)^3}
{\eta_{100,15}(\tau)^2\eta_{100,35}(\tau)^2\eta_{100,40}(\tau)\eta_{100,50}(\tau)}
+\frac{8\eta_{100, 5}(\tau)^4\eta_{100, 20}(\tau)^4\eta_{100, 25}(\tau)^2\eta_{100, 30}(\tau)^4\eta_{100, 45}(\tau)^4}
{\eta_{100, 10}(\tau)^4\eta_{100, 15}(\tau)\eta_{100, 35}(\tau)\eta_{100, 40}(\tau)^4}
\\[5pt]
&-\frac{16\eta_{10,5}(\tau)^{\frac{1}{2}}\eta_{100,15}(\tau)^2\eta_{100,35}(\tau)^2}{\eta_{100,25}(\tau)}
+\frac{24\eta_{10,5}(\tau)^{\frac{1}{2}}\eta_{100,15}(\tau)^2\eta_{100,20}(\tau)\eta_{100,30}(\tau)^2\eta_{100,35}(\tau)^2}
{\eta_{100,10}(\tau)\eta_{100,25}(\tau)\eta_{100,40}(\tau)\eta_{100,50}(\tau)}
\\[5pt]
&+\frac{8\eta_{10,0}(\tau)\eta_{100,5}(\tau)\eta_{100,10}(\tau)\eta_{100,40}(\tau)\eta_{100,45}(\tau)}
{\eta_{100,0}(\tau)\eta_{100,15}(\tau)\eta_{100,20}(\tau)\eta_{100,35}(\tau)\eta_{100,50}(\tau)^2}
-\frac{24\eta_{100, 5}(\tau)^2\eta_{100, 20}(\tau)^2\eta_{100, 25}(\tau)^4\eta_{100, 30}(\tau)^2\eta_{100, 45}(\tau)^2}
{\eta_{100, 10}(\tau)\eta_{100, 40}(\tau)^2\eta_{100, 50}(\tau)}
\\[5pt]
&-\frac{32\eta_{100, 5}(\tau)^3\eta_{100, 20}(\tau)\eta_{100, 25}(\tau)^2\eta_{100, 30}(\tau)^2\eta_{100, 45}(\tau)^3}
{\eta_{100, 10}(\tau)\eta_{100, 40}(\tau)\eta_{100, 50}(\tau)}
-\frac{16\eta_{100, 5}(\tau)\eta_{100, 10}(\tau)^3\eta_{100, 20}(\tau)^2\eta_{100, 30}(\tau)^3\eta_{100, 45}(\tau)}
{\eta_{100, 15}(\tau)\eta_{100, 35}(\tau)}
\\[5pt]
&+\frac{2\eta_{10,5}(\tau)^{\frac{3}{2}}\eta_{100,20}(\tau)^2\eta_{100,25}(\tau)^5\eta_{100,30}(\tau)^3}
{\eta_{100,15}(\tau)^2\eta_{100,35}(\tau)^2\eta_{100,50}(\tau)}
-\frac{16\eta_{10,5}(\tau)\eta_{100,10}(\tau)\eta_{100,40}(\tau)}
{\eta_{100,5}(\tau)\eta_{100,20}(\tau)\eta_{100,30}(\tau)\eta_{100,45}(\tau)}
\\[5pt]
&+\frac{24\eta_{100, 5}(\tau)^2\eta_{100, 20}(\tau)\eta_{100, 25}(\tau)^4\eta_{100, 30}(\tau)^2\eta_{100, 45}(\tau)^2}
{\eta_{100, 10}(\tau)\eta_{100, 40}(\tau)\eta_{100, 50}(\tau)}
+\frac{8\eta_{100, 5}(\tau)^6\eta_{100, 20}(\tau)^4\eta_{100, 25}(\tau)^2\eta_{100, 30}(\tau)^4\eta_{100, 45}(\tau)^6}
{\eta_{100, 10}(\tau)^4\eta_{100, 15}(\tau)^3\eta_{100, 35}(\tau)^3\eta_{100, 40}(\tau)^4}
\\[5pt]
&-\frac{32\eta_{10,5}(\tau)\eta_{100,20}(\tau)^3\eta_{100,30}(\tau)^2}
{\eta_{100,10}(\tau)\eta_{100,25}(\tau)^2\eta_{100,40}(\tau)^3\eta_{100,50}(\tau)}
-\frac{32\eta_{100, 5}(\tau)^5\eta_{100, 20}(\tau)^2\eta_{100, 30}(\tau)^2\eta_{100, 45}(\tau)^5}
{\eta_{100, 10}(\tau)\eta_{100, 15}(\tau)\eta_{100, 35}(\tau)\eta_{100, 40}(\tau)^2\eta_{100, 50}(\tau)}
\\[5pt]
&+\frac{4\eta_{100, 5}(\tau)\eta_{100, 20}(\tau)\eta_{100, 25}(\tau)^6\eta_{100, 30}(\tau)^2\eta_{100, 45}(\tau)}
{\eta_{100, 10}(\tau)\eta_{100, 40}(\tau)\eta_{100, 50}(\tau)}
-\frac{4\eta_{100, 10}(\tau)^3\eta_{100, 20}(\tau)^2\eta_{100, 25}(\tau)^4\eta_{100, 30}(\tau)^3}
{\eta_{100, 5}(\tau)\eta_{100, 15}(\tau)\eta_{100, 35}(\tau)\eta_{100, 45}(\tau)}
\\[5pt]
&+\frac{32\eta_{100, 5}(\tau)^4\eta_{100, 20}(\tau)\eta_{100, 30}(\tau)\eta_{100, 45}(\tau)^4}
{\eta_{100, 10}(\tau)\eta_{100, 40}(\tau)}
-\frac{32\eta_{10,5}(\tau)\eta_{100,20}(\tau)^2\eta_{100,30}(\tau)^2}
{\eta_{100,10}(\tau)\eta_{100,25}(\tau)^2\eta_{100,40}(\tau)^2\eta_{100,50}(\tau)}
\\[5pt]
&+\frac{16\eta_{100, 5}(\tau)^4\eta_{100, 20}(\tau)^2\eta_{100, 25}(\tau)^2\eta_{100, 30}(\tau)^2\eta_{100, 45}(\tau)^4}
{\eta_{100, 10}(\tau)\eta_{100, 15}(\tau)\eta_{100, 35}(\tau)\eta_{100, 40}(\tau)^2\eta_{100, 50}(\tau)}
+\frac{32\eta_{100, 5}(\tau)^5\eta_{100, 20}(\tau)\eta_{100, 30}(\tau)^2\eta_{100, 45}(\tau)^5}
{\eta_{100, 10}(\tau)\eta_{100, 15}(\tau)\eta_{100, 35}(\tau)\eta_{100, 40}(\tau)\eta_{100, 50}(\tau)}
\\[5pt]
&+\frac{64\eta_{10,5}(\tau)\eta_{100,10}(\tau)\eta_{100,20}(\tau)\eta_{100,30}(\tau)}
{\eta_{100,25}(\tau)^2\eta_{100,40}(\tau)\eta_{100,50}(\tau)}\Bigg).
\end{align*}
\end{footnotesize}
In light of Theorem \ref{Robins-1994}, it can be verified that each term of the above identity is a modular function with respect to $\Gamma_1(100)$. Employing the algorithm  by Garvan and Liang \cite{Garvan-Liang-2016}, we deduce that the constant $B$  in \eqref{constantB} is equal to  $-640$. By Theorem \ref{eta-lem},
\eqref{generfun-2} can be justified if
it holds for coefficients of $q^{n}$ for $1\leq n \leq 640$,
which can be easily checked. This completes the proof. \qed

%Theorem  \ref{main-result-10-1-M2} can be deduced from Lemma \ref{U1-2q^3U2}.

 Theorems \ref{main-result-10-1-M2} and  \ref{main-result-10-2-M2}
 are immediate consequences of Lemmas \ref{U1-2q^3U2} and \ref{V1-U1+2q^3U2-2q^7V2}.

\vskip 0.2cm
\noindent{\it Proof of Theorem  \ref{main-result-10-1-M2}.}
The relation \eqref{main-result-10-M2-F1} in Lemma \ref{main-result-10-lem-1} states that
\begin{align}\label{main-result-10-1-pff}
&\sum_{n=0}^\infty(\overline{N}_2(0,10,n)+\overline{N}_2(1,10,n)
-\overline{N}_2(4,10,n)-\overline{N}_2(5,10,n))q^n
\nonumber\\
&\quad =\frac{2(-q;q)_\infty}{(q;q)_\infty}\left(\frac{1}{2}U_1-q^2U_3\right)+
\frac{2q^{24}}{J_{25,50}}
\sum_{n=-\infty}^\infty\frac{(-1)^nq^{25n^2+50n}}{1+q^{50n+40}}.
\end{align}
Subsituting \eqref{U1-2q^3U2-eq} into \eqref{main-result-10-1-pff}, we obtain that
\begin{align*}
&\sum_{n=0}^\infty(\overline{N}_2(0,10,n)+\overline{N}_2(1,10,n)-\overline{N}_2(4,10,n)-\overline{N}_2(5,10,n))q^n
\nonumber\\[5pt]
&=\frac{2J_{50}^4J_{25,50}}
{J_{5,10}^6}\left(A_0(q^5)+qA_1(q^5)+q^2A_2(q^5)+q^3A_3(q^5)\right)
\nonumber\\[5pt]
&\quad\quad\quad+2q^4\left(\frac{J_{50}^4J_{25,50}}
{J_{5,10}^6}A_4(q^5)+\frac{q^{20}}{J_{25,50}}
\sum_{n=-\infty}^\infty\frac{(-1)^nq^{25n^2+50n}}{1+q^{50n+40}}\right).
\end{align*}
Comparing the coefficients $q^{5n}$, $q^{5n+1}$, $q^{5n+2}$, $q^{5n+3}$ and $q^{5n+4}$ of above relation, we get \eqref{main-result-eq1-10-M2-1} and \eqref{main-result-eq1-10-M2-2}. This completes the proof. \qed
\vskip 0.2cm

%We now present a proof of Theorem  \ref{main-result-10-2-M2} by using Lemma \ref{V1-U1+2q^3U2-2q^7V2}.

\noindent{\it Proof of Theorem  \ref{main-result-10-2-M2}.}
The relation  \eqref{main-result-10-M2-F2} in Lemma \ref{main-result-10-lem-1} reads
\begin{align}\label{main-result-10-2-pfeq}
&\sum_{n=0}^\infty(\overline{N}_2(1,10,n)+\overline{N}_2(2,10,n)
-\overline{N}_2(3,10,n)-\overline{N}_2(4,10,n))q^n\nonumber \\[5pt]
&\quad=\frac{2(-q;q)_\infty}{(q;q)_\infty}\left(\frac{1}{2}V_1-q^7V_2\right)-
\frac{2q^{21}}{J_{25,50}}
\sum_{n=-\infty}^\infty\frac{(-1)^nq^{25n^2+50n}}{1+q^{50n+30}}
\nonumber \\[5pt]
&\quad\quad\quad -\left\{\frac{2(-q;q)_\infty}{(q;q)_\infty}\left(\frac{1}{2}U_1-q^3U_2\right)+
\frac{2q^{24}}{J_{25,50}}
\sum_{n=-\infty}^\infty\frac{(-1)^nq^{25n^2+50n}}{1+q^{50n+40}}\right\}
\nonumber\\[5pt]
&\quad=\frac{2(-q;q)_\infty}{(q;q)_\infty}\left(\frac{1}{2}V_1-\frac{1}{2}U_1+q^3U_2-q^7V_2\right)
\nonumber\\[5pt]
&\quad\quad\quad -
\frac{2q^{21}}{J_{25,50}}
\sum_{n=-\infty}^\infty\frac{(-1)^nq^{25n^2+50n}}{1+q^{50n+30}}
-
\frac{2q^{24}}{J_{25,50}}
\sum_{n=-\infty}^\infty\frac{(-1)^nq^{25n^2+50n}}{1+q^{50n+40}}.
\end{align}
Substituting \eqref{V1-U1+2q^3U2-2q^7V2-eq} into \eqref{main-result-10-2-pfeq}, we obtain that
\begin{align*}
&\sum_{n=0}^\infty(\overline{N}_2(1,10,n)+\overline{N}_2(2,10,n)
-\overline{N}_2(3,10,n)-\overline{N}_2(4,10,n))q^n
\nonumber\\
&\quad=\frac{2J_{50}^4J_{25,50}}
{J_{5,10}^6}B_0(q^5)
+2q\left(\frac{J_{50}^4J_{25,50}}
{J_{5,10}^6}B_1(q^5)
-\frac{q^{20}}{J_{25,50}}
\sum_{n=-\infty}^\infty\frac{(-1)^nq^{25n^2+50n}}{1+q^{50n+30}}\right)
\nonumber\\
&\quad\quad\quad
+2q^2\frac{J_{50}^4J_{25,50}}
{J_{5,10}^6}B_2(q^5)
+2q^3\frac{J_{50}^4J_{25,50}}
{J_{5,10}^6}B_3(q^5)
\nonumber\\
&\quad\quad\quad
+2q^4\left(\frac{J_{50}^4J_{25,50}}
{J_{5,10}^6}B_4(q^5)-\frac{q^{20}}{J_{25,50}}
\sum_{n=-\infty}^\infty\frac{(-1)^nq^{25n^2+50n}}{1+q^{50n+40}}\right).
\end{align*}
Comparing the coefficients $q^{5n}$, $q^{5n+1}$, $q^{5n+2}$, $q^{5n+3}$ and $q^{5n+4}$ of above relation, we get \eqref{main-result-eq2-10-M2-1}, \eqref{main-result-eq2-10-M2-2} and \eqref{main-result-eq2-10-M2-3}. This completes the proof. \qed

\section{Proof of Theorem \ref{rank-ineq1}}

Before presenting the proof of Theorem \ref{rank-ineq1}, we recall the following theorem due to Liaw \cite{Liaw-1999}.

\begin{thm}\label{posi}
If $p$ and $r$ are positive integers with $p\geq 2$ and $r<p$ and
\[\sum_{n=0}^\infty b_{p,r}(n)q^n=\frac{(q^p;q^p)_\infty}{(q^r;q^p)_\infty(q^{p-r};q^p)_\infty},\]
then $b_{p,r}(n)\geq 0$ for all $n$.
\end{thm}

We are now in a position to prove Theorem \ref{rank-ineq1}. \vskip 0.2cm

\noindent{\it Proof of Theorem \ref{rank-ineq1}.} We first prove \eqref{ineq-3-0-1-f}, \eqref{ineq-3-0-3-f} and  \eqref{ineq-3-0-2-f}.
From Theorem \ref{main-result-6-M2}, we find that
\begin{align}\label{neq-3-0-1-pf-t}
\overline{R}_{0,2}(0,6)
=\frac{J_{6}^9J_{2,4}^2}{J^6_{1,6}J_{2}^3J_{3,6}}
=\frac{(q^6;q^6)_\infty^2}{(q,q^5;q^6)_\infty^2}
\frac{(-q;q^2)_\infty}{(q^4;q^4)_\infty(q^3;q^6)_\infty(q,q^5;q^6)_\infty^3}.
\end{align}
By Theorem \ref{posi}, we see that the coefficients of $q^n$ in
\[\frac{(q^6;q^6)_\infty^2}{(q,q^5;q^6)_\infty^2}\]
is nonnegative for $n\geq 0$. This implies that
\begin{align}\label{ineq-3-0-1}
\overline{N}_2(0,6,3n)+\overline{N}_2(1,6,3n)
  &\geq\overline{N}_2(2,6,3n)+\overline{N}_2(3,6,3n).
\end{align}
Since $i\equiv s\pmod{m}$ if and only if  $i\equiv s\pmod{2m}$ or $i\equiv m+s\pmod{2m}$, we find that
\begin{equation}\label{rank-sym}
\overline{N}_2(s,m,n)=\overline{N}_2(s,2m,n)+\overline{N}_2(m+s,2m,n),
\end{equation}
and thus
\begin{equation}\label{N2-1}
\overline{N}_2(0,3,3n)=\overline{N}_2(0,6,3n)+\overline{N}_2(3,6,3n).
\end{equation}
Moreover, by \eqref{Sy-N}
\[ \overline{N}_2(m+s,2m,n)=\overline{N}_2(m-s,2m,n),\]
\eqref{rank-sym} can be rewritten as
\begin{equation}\label{rank-sym-1}
\overline{N}_2(s,m,n)
=\overline{N}_2(s,2m,n)+\overline{N}_2(m-s,2m,n).
\end{equation}
In particular, we get
\begin{equation}\label{N2-2}
\overline{N}_2(1,3,3n)=\overline{N}_2(1,6,3n)+\overline{N}_2(2,6,3n).
\end{equation}
 By means of the following identity due to  Lovejoy and Osburn \cite[Theorem 1.1]{Lovejoy-Osburn-2010},
\begin{align*}
  &\sum_{n=0}^\infty(\overline{N}_2(0,3,3n)-\overline{N}_2(1,3,3n))q^n=-1+\frac{J_{6}^3J_{3,6}}{J^2_{1,6}J_{2}},
  %\nonumber\\
  %&=\sum_{n=0}^\infty(\overline{N}_2(0,6,3n)+\overline{N}_2(3,6,3n)-\overline{N}_2(1,6,3n)-\overline{N}_2(2,6,3n))q^n
  %\nonumber\\
\end{align*}
together with \eqref{N2-1} and \eqref{N2-2}, we obtain
\begin{align}\label{Lovejoy-Osburn-2010-3-0}
\sum_{n=0}^\infty(\overline{N}_2(0,6,3n)+\overline{N}_2(3,6,3n)-\overline{N}_2(1,6,3n)-\overline{N}_2(2,6,3n))q^n
=-1+\frac{J_{6}^3J_{3,6}}{J^2_{1,6}J_{2}}.
\end{align}
Consequently,
\begin{align*}
 \frac{J_{6}^3J_{3,6}}{J^2_{1,6}J_{2}}
=\frac{(q^3;q^6)_\infty^2(q^6;q^6)_\infty}{(q,q^5;q^6)_\infty^2(q^2,q^4;q^6)_\infty}
=\frac{(q^3;q^3)_\infty}{(q,q^2;q^3)_\infty}\frac{(q^3;q^6)_\infty}{(q,q^5;q^6)_\infty}
=\frac{(q^3;q^3)_\infty}{(q,q^2;q^3)_\infty}\frac{(q^3;q^6)^2_\infty}{(q;q^2)_\infty}
\end{align*}
and
\[\frac{(q^3;q^6)^2_\infty}{(q;q^2)_\infty}=\frac{(-q;q)_\infty}{(-q^3;q^3)_\infty^2}
=\frac{(-q;q^3)_\infty(-q^2;q^3)_\infty}{(-q^3;q^3)_\infty}
%=\frac{(-q;q^3)_\infty(-q^2;q^3)_\infty(q^3;q^3)_\infty}{(-q^3;q^3)_\infty(q^3;q^3)_\infty}
=\frac{1}{(q^6;q^6)_\infty}\sum_{n=-\infty}^\infty q^{(3n^2-n)/2}.\]
Therefore,  \eqref{Lovejoy-Osburn-2010-3-0} can be written as
\begin{align*}
&\sum_{n=0}^\infty(\overline{N}_2(0,6,3n)+\overline{N}_2(3,6,3n)
-\overline{N}_2(1,6,3n)-\overline{N}_2(2,6,3n))q^n
\\&=-1+\frac{(q^3;q^3)_\infty}{(q,q^2;q^3)_\infty}\frac{1}{(q^6;q^6)_\infty}\sum_{n=-\infty}^\infty q^{(3n^2-n)/2}.
\end{align*}
By Theorem \ref{posi}, we deduce that for $n\geq1$
\begin{align}\label{ineq-3-0-2-1}
\overline{N}_2(0,6,3n)+\overline{N}_2(3,6,3n)\geq\overline{N}_2(1,6,3n)+\overline{N}_2(2,6,3n),
\end{align}
so \eqref{ineq-3-0-2-f} is verified. Adding \eqref{ineq-3-0-1} to \eqref{ineq-3-0-2-1}, we reach \eqref{ineq-3-0-1-f}.

Subtracting \eqref{Lovejoy-Osburn-2010-3-0} from \eqref{neq-3-0-1-pf-t}, we obtain that
\begin{align*}
\sum_{n=0}^\infty \left(\overline{N}_2(1,6,3n)-\overline{N}_2(3,6,3n)\right)q^n
&=\frac{1}{2}\cdot\frac{J_{6}^9J_{2,4}^2}{J^6_{1,6}J_{2}^3J_{3,6}}+\frac{1}{2}-\frac{1}{2}\cdot\frac{J_{6}^3J_{3,6}}{J^2_{1,6}J_{2}}
\\[3pt]
&=\frac{1}{2}+
\frac{1}{2}\cdot\frac{J_{6}^3J_{3,6}}{J^2_{1,6}J_{2}}\left(\frac{J_{6}^6J_{2,4}^2}{J^4_{1,6}J_{2}^2J_{3,6}^2}-1\right)
\\[3pt]
&=\frac{1}{2}+
\frac{1}{2}\cdot\frac{J_{6}^3J_{3,6}}{J^2_{1,6}J_{2}}\cdot\left((-q;q^2)_\infty^2(-q;q)_\infty^2-1\right),
\end{align*}
which yields \eqref{ineq-3-0-3-f}.

To verify \eqref{ineq-3-1-1-f},  we rewrite \eqref{main-result-6-M2-2} as
\begin{align}\label{neq-3-1-1-pf-t}
  &\overline{R}_{0,2}(1,6)=\frac{2J_{6}^9J_{2,4}^2}{J^5_{1,6}J_{2}^3J_{3,6}^2}
 =\frac{2(q^6;q^6)_\infty^2}{(q,q^5;q^6)_\infty^2}
 \frac{(-q;q^2)_\infty}{(q^4;q^4)_\infty(q,q^5;q^6)_\infty^2(q^3;q^6)_\infty^3}.
\end{align}
By Theorem \ref{posi}, we see that the coefficients of $q^n$ in
\[\frac{(q^6;q^6)_\infty^2}{(q,q^5;q^6)_\infty^2}\]
is nonnegative for $n\geq 0$, which implies that
\begin{align}\label{ineq-3-1-1}
\overline{N}_2(0,6,3n+1)+\overline{N}_2(1,6,3n+1)
  &\geq\overline{N}_2(2,6,3n+1)+\overline{N}_2(3,6,3n+1).
\end{align}
Using the following identity due to Lovejoy and Osburn \cite[Theorem 1.1]{Lovejoy-Osburn-2010}:
\begin{align*}
&\sum_{n=0}^\infty(\overline{N}_2(0,3,3n+1)-\overline{N}_2(1,3,3n+1))q^n
=\frac{2J_{6}^3}{J_{1,6}J_{2}}
=\frac{(q^6;q^6)_\infty}{(q^2,q^4;q^6)_\infty}\frac{2}{(q,q^5;q^6)_\infty},
\end{align*}
together with the relations \eqref{rank-sym} and \eqref{rank-sym-1}, we get
\begin{align}\label{Lovejoy-Osburn-2010-3-1}
&\sum_{n=0}^\infty(\overline{N}_2(0,6,3n+1)+\overline{N}_2(3,6,3n+1)
-\overline{N}_2(1,6,3n+1)-\overline{N}_2(2,6,3n+1))q^n
\nonumber\\
&\quad=\frac{(q^6;q^6)_\infty}{(q^2,q^4;q^6)_\infty}\frac{2}{(q,q^5;q^6)_\infty}.
\end{align}
By Theorem \ref{posi}, we see that the coefficients of $q^n$ in
\[\frac{(q^6;q^6)_\infty}{(q^2,q^4;q^6)_\infty}\]
is nonnegative for $n\geq 0$.
Hence
\begin{align}\label{ineq-3-1-2-1}
\overline{N}_2(0,6,3n+1)+\overline{N}_2(3,6,3n+1)\geq\overline{N}_2(1,6,3n+1)+\overline{N}_2(2,6,3n+1),
\end{align}
as given in \eqref{ineq-3-1-2-f}.
 Adding \eqref{ineq-3-1-1} to \eqref{ineq-3-1-2-f},  we arrive at \eqref{ineq-3-1-1-f}.
Subtracting \eqref{Lovejoy-Osburn-2010-3-1} from \eqref{neq-3-1-1-pf-t}, we obtain that
\begin{align}\label{ineq-3---1}
\sum_{n=0}^\infty \left(\overline{N}_2(1,6,3n+1)-\overline{N}_2(3,6,3n+1)\right)q^n
&=\frac{J_{6}^9J_{2,4}^2}{J^5_{1,6}J_{2}^3J_{3,6}^2}-\frac{J_{6}^3}{J_{1,6}J_{2}}
\nonumber\\[3pt]
&=
\frac{J_{6}^3}{J_{1,6}J_{2}}\left(\frac{J_{6}^6J_{2,4}^2}{J^4_{1,6}J_{2}^2J_{3,6}^2}-1\right)
\nonumber\\[3pt]
&=\frac{J_{6}^3}{J_{1,6}J_{2}}\cdot\left((-q;q^2)_\infty^2(-q;q)_\infty^2-1\right).
\end{align}
It is easily seen that the coefficients of $q^n$ in
\eqref{ineq-3---1} are nonnegative for $n\geq1$, and this justifies \eqref{ineq-3-1-3-f}.

It remains to prove \eqref{ineq-3-2-1-f}. Using the following
relation due to Lovejoy and Osburn \cite[Theorem 1.1]{Lovejoy-Osburn-2010},
\begin{align*}
&\sum_{n=0}^\infty(\overline{N}_2(0,3,3n+2)-\overline{N}_2(1,3,3n+2))q^n
\\&\quad=\frac{4(q^6;q^6)_\infty^4}{(q^2;q^2)_\infty(q^3;q^3)_\infty^2}
+\frac{6(-q^3;q^3)_\infty}{(q^3;q^3)_\infty}\sum_{n=-\infty}^\infty\frac{(-1)^n q^{3n^2+6n+1}}{1-q^{6n+2}}
\\&\quad=\frac{(q^6;q^6)_\infty}{(q^2,q^4;q^6)_\infty}\frac{4}{(q^3;q^6)_\infty^2}
+\frac{6(-q^3;q^3)_\infty}{(q^3;q^3)_\infty}\sum_{n=-\infty}^\infty\frac{(-1)^n q^{3n^2+6n+1}}{1-q^{6n+2}},
\end{align*}
together with \eqref{N2-1} and \eqref{N2-2}, we get
\begin{align*}
&\sum_{n=0}^\infty(\overline{N}_2(0,6,3n+2)+\overline{N}_2(3,6,3n+2)-\overline{N}_2(1,6,3n+2)-\overline{N}_2(2,6,3n+2))q^n
\\
&\quad=\frac{(q^6;q^6)_\infty}{(q^2,q^4;q^6)_\infty}\frac{4}{(q^3;q^6)_\infty^2}
+\frac{6(-q^3;q^3)_\infty}{(q^3;q^3)_\infty}\sum_{n=-\infty}^\infty\frac{(-1)^n q^{3n^2+6n+1}}{1-q^{6n+2}}.
\end{align*}
By Theorem \ref{posi}, we see that the coefficient  of $q^n$ in
\[\frac{(q^6;q^6)_\infty}{(q^2,q^4;q^6)_\infty}\]
is nonnegative for $n\geq 0$.

We now consider the positivity of
\begin{align}\label{ineq-3---2}
\frac{(-q^3;q^3)_\infty}{(q^3;q^3)_\infty}\sum_{n=-\infty}^\infty\frac{(-1)^n q^{3n^2+6n+1}}{1-q^{6n+2}}.
\end{align}
The sum in \eqref{ineq-3---2} equals
\begin{align}\label{ineq-3---2-1}
&\sum_{n=0}^\infty\frac{(-1)^nq^{3n^2+6n+1}}{1-q^{6n+2}}
+\sum_{n=0}^\infty\frac{(-1)^nq^{3n^2+6n+2}}{1-q^{6n+4}}
\nonumber\\[3pt]
&\quad=\sum_{n=0}^\infty\frac{q^{12n^2+12n+1}}{1-q^{12n+2}}
-\sum_{n=0}^\infty\frac{q^{12n^2+24n+10}}{1-q^{12n+8}}
+\sum_{n=0}^\infty\frac{q^{12n^2+12n+2}}{1-q^{12n+4}}
-\sum_{n=0}^\infty\frac{q^{12n^2+24n+11}}{1-q^{12n+10}}
\nonumber\\[3pt]
&\quad=\sum_{n=0}^\infty\frac{q^{12n^2+12n+1}}{(1-q^{12n+2})(1-q^{12n+8})}
\left((1-q^{12n+8})(1-q^{12n+9})+q^{24n+11}(1-q^6)\right)
\nonumber\\[3pt]
&\hskip1.5cm+\sum_{n=0}^\infty\frac{q^{12n^2+12n+2}}{(1-q^{12n+4})(1-q^{12n+10})}
\left((1-q^{12n+9 })(1-q^{12n+10})+q^{24n+13}(1-q^6)\right)
\nonumber\\[3pt]
&\quad=(1-q^3)\left(\sum_{n=0}^\infty\frac{q^{12n^2+12n+1}}{1-q^{12n+2}}\frac{1-q^{12n+9 }}{1-q^3}+\sum_{n=0}^\infty\frac{q^{12n^2+36n+12}(1+q^3)}{(1-q^{12n+2})(1-q^{12n+8})}\right.
\nonumber\\[3pt]
&\hskip3cm\left.+\sum_{n=0}^\infty\frac{q^{12n^2+12n+2}}{1-q^{12n+4}}\frac{1-q^{12n+9 }}{1-q^3}+\sum_{n=0}^\infty\frac{q^{12n^2+36n+15}(1+q^3)}{(1-q^{12n+4})(1-q^{12n+10})}\right)
\nonumber\\[3pt]
&\quad=(1-q^3)\left(\sum_{n=0}^\infty\frac{q^{12n^2+12n+1}}{1-q^{12n+2}}
\sum_{m=0}^{4n+2}q^{3m}+\sum_{n=0}^\infty\frac{q^{12n^2+36n+12}(1+q^3)}
{(1-q^{12n+2})(1-q^{12n+8})}\right.
\nonumber\\[3pt]
&\hskip3cm\left.+\sum_{n=0}^\infty
\frac{q^{12n^2+12n+2}}{1-q^{12n+4}}\sum_{m=0}^{4n+2}q^{3m}
+\sum_{n=0}^\infty\frac{q^{12n^2+36n+15}(1+q^3)}{(1-q^{12n+4})(1-q^{12n+10})}\right).
\end{align}
Noting that
\begin{align}\label{ineq-3---2-2}
\frac{(-q^3;q^3)_\infty}{(q^3;q^3)_\infty}(1-q^3)
=\frac{(-q^3;q^3)_\infty}{(q^6;q^3)_\infty},
\end{align}
it follows from \eqref{ineq-3---2-1} and \eqref{ineq-3---2-2} that
the coefficient of $q^n$ in \eqref{ineq-3---2} is nonnegative for $n\geq 0$.
This completes the proof. \qed

\section{Proofs of Theorems \ref{mock-1} and \ref{mock-2}}
This section is devoted to the proofs of the relations between the $M_2$-rank differences of overpartitions and mock theta functions as stated in Theorems \ref{mock-1} and \ref{mock-2}.
It is known that mock theta  functions can be expressed in terms of  the Appell-Lerch sum $m(x,q,z)$.
Recall that the Appell-Lerch sum is defined by
\begin{equation}
m(x,q,z)=\frac{1}{j(z;q)}\sum_{r=-\infty }^{\infty}
\frac{(-1)^r q^{r\choose 2}z^r}{1-q^{r-1}xz},
\end{equation}
where $x,z \in \mathbb{C}^*$ with neither $z$ nor $xz$ an integral power of $q$.

Hickerson and  Mortenson \cite{Hickerson-Mortenson-2014} proved that
\begin{align}
%\label{a2}
%\omega(q)&=&-2q^{-1}m(q,q^{6},q^2)+\frac{J^3_{6}}{J_{2}J_{3,6}},\\[3pt]
\label{a1}
\rho(q)&=q^{-1}m(q,q^{6},-q),\\[3pt]\label{a3}
\phi(q)&=-q^{-1}m(q,q^{10},q)-q^{-1}m(q,q^{10},q^2),\\[3pt]\label{a4}
\psi(q)&=-m(q^3,q^{10},q)-m(q^3,q^{10},q^3).
\end{align}

We are ready to present the proofs of Theorems \ref{mock-1} and \ref{mock-2}.  \vskip 0.2cm

\noindent{\it Proof of Theorem \ref{mock-1}.}
From Theorem \ref{main-result-6-M2}, we have
\begin{align}
\overline{R}_{0,2}(2,6)
&=\frac{4J_{6}^9J_{2,4}^2}{J^4_{1,6}J^3_{2}J^3_{3,6}}
+\frac{2q^2}{J_{3,6}}
\sum_{n=-\infty}^\infty\frac{(-1)^nq^{3n^2+6n}}{1+q^{6n+4}}\\[3pt]\nonumber
&=\frac{4J_{6}^9J_{2,4}^2}{J^4_{1,6}J^3_{2}J^3_{3,6}}
-\frac{2q^2}{J_{3,6}}
\sum_{n=-\infty}^\infty\frac{(-1)^nq^{3n^2-3}}{1+q^{6n-2}}
\\[3pt]\label{main-result-6aaaa}
&=\frac{4J_{6}^9J_{2,4}^2}{J^4_{1,6}J^3_{2}J^3_{3,6}}
-2q^{-1}m(-q,q^{6},q^3).
\end{align}
Replacing $q$ by $-q$ in \eqref{a1}, we deduce that
\begin{align}\label{a-1}
\rho(-q)=-q^{-1}m(-q,q^{6},q).
\end{align}
Recall the following identity in Hickerson and  Mortenson  \cite{Hickerson-Mortenson-2014}:
\begin{equation}\label{mm}
m(x,q,z_1)-m(x,q,z_0)=\frac{z_0J_1^3j(z_1/z_0;q)j(xz_0z_1;q)}
{j(z_0;q)j(z_1;q)j(xz_0;q)j(xz_1;q)}.
\end{equation}
Replacing $q$, $z_0$, $z_1$ and $x$  with $q^6$, $q$, $q^3$ and $-q$, respectively, in \eqref{mm}, we get
\begin{equation}\label{m3}
m(-q,q^6,q^3)-m(-q,q^6,q)=\frac{q J^3_{6}J_2\overline{J}_{1,6}}{J_{1,6}J_{3,6}\overline{J}_{2}^2}.
\end{equation}
Combining \eqref{a-1} and \eqref{m3}, we obtain
\begin{align}\label{mock--1}
m(-q,q^6,q^3)=-q\rho(-q)+\frac{q J^3_{6}J_2\overline{J}_{1,6}}{J_{1,6}J_{3,6}\overline{J}_{2}^2}.
\end{align}
Substituting \eqref{mock--1} into \eqref{main-result-6aaaa}, we get
\eqref{mock-1-eqn}.
This completes the proof of Theorem \ref{mock-1}.\qed

We conclude with a proof of Theorem \ref{mock-2}. \vskip 0.2cm

\noindent{\it Proof of Theorem \ref{mock-2}.}
To prove \eqref{mock-2-1}, we apply Theorem \ref{main-result-10-1-M2} to get
\begin{align}\label{main-result-10-1-tempt}
\overline{R}_{0,4}(4,10)
&=\frac{2J_{10}^4J_{5,10}}{J_{1,2}^6}A_4(q)
+\frac{2q^4}{J_{5,10}}\sum_{n=-\infty}^\infty\frac{(-1)^n q^{5n^2+10n}}{1+q^{10n+8}}
\nonumber\\[3pt]
&=\frac{2J_{10}^4J_{5,10}}{J_{1,2}^6}A_4(q)
-\frac{2}{qJ_{5,10}}\sum_{n=-\infty}^\infty\frac{(-1)^n q^{5n^2}}{1+q^{10n-2}}
\nonumber\\[3pt]
&=\frac{2J_{10}^4J_{5,10}}{J_{1,2}^6}A_4(q)
-2q^{-1}m(-q^3,q^{10},q^5),
\end{align}
where $A_4(q)$ is defined as \eqref{10-eq-A4}.
Replacing $q$ by $-q$ in \eqref{a4}, we are led to
\begin{align}\label{psi-q}
\psi(-q)=-m(-q^3,q^{10},-q)-m(-q^3,q^{10},-q^3).
\end{align}
Replacing $q$, $z_0$, $z_1$ and $x$ by $q^{10}$, $-q$, $q^5$ and $-q^3$ in \eqref{mm}, respectively, we find that
\begin{equation}\label{mm-1}
  m(-q^3,q^{10},q^5)-m(-q^3,q^{10},-q)=-\frac{qJ_{10}^3\overline{J}_{4,10}J_{1,10}}
{J_{4,10}J_{5,10}\overline{J}_{1,10}\overline{J}_{2,10}}.
\end{equation}
Replacing $q$, $z_1$, $z_0$ and $x$ in \eqref{mm} by $q^{10}$, $q^{5}$, $-q^{3}$ and $-q^{3}$, respectively, we see that
\begin{equation}\label{mm-2}
  m(-q^3,q^{10},q^5)-m(-q^3,q^{10},-q^3)=\frac{q^2J_{10}^3J_{1,10}\overline{J}_{2,10}}
  {J_{4,10}J_{5,10}\overline{J}_{2,10}\overline{J}_{3,10}}.
\end{equation}
Combining \eqref{mm-1} and \eqref{mm-2} gives
\begin{align}\label{mock-1---3}
m(-q^3,q^{10},q^5)&=\frac{1}{2}m(-q^3,q^{10},-q)
+\frac{1}{2}m(-q^3,q^{10},-q^3)
\nonumber\\[3pt]
&\quad+\frac{q^2J_{10}^3J_{1,10}\overline{J}_{2,10}}
  {2J_{4,10}J_{5,10}\overline{J}_{2,10}\overline{J}_{3,10}}
  -\frac{qJ_{10}^3\overline{J}_{4,10}J_{1,10}}
{2J_{4,10}J_{5,10}\overline{J}_{1,10}\overline{J}_{2,10}}.
\end{align}
 Utilizing
 \eqref{psi-q} and \eqref{mock-1---3}, $m(-q^3,q^{10},q^5)$ can be written as
\begin{align}\label{mock-1---2}
m(-q^3,q^{10},q^5)=-\frac{1}{2}\psi(-q)
+\frac{q^2J_{10}^3J_{1,10}\overline{J}_{2,10}}
  {2J_{4,10}J_{5,10}\overline{J}_{2,10}\overline{J}_{3,10}}
  -\frac{qJ_{10}^3\overline{J}_{4,10}J_{1,10}}
{2J_{4,10}J_{5,10}\overline{J}_{1,10}\overline{J}_{2,10}}.
\end{align}
Substituting \eqref{mock-1---2} into \eqref{main-result-10-1-tempt}, we implies
\eqref{mock-2-1}.

We now consider \eqref{mock-2-2}. By Theorem \ref{main-result-10-2-M2},
we find that
\begin{align}\label{main-result-10-2-tempt}
\overline{R}_{1,3}(1,10)
&=\frac{2J_{10}^4J_{5,10}}{J_{1,2}^6}B_1(q)
-\frac{2q^4}{J_{5,10}}\sum_{n=-\infty}^\infty\frac{(-1)^n q^{5n^2+10n}}{1+q^{10n+6}}
\nonumber\\[3pt]
&=\frac{2J_{10}^4J_{5,10}}{J_{1,2}^6}B_1(q)
+\frac{2}{qJ_{5,10}}\sum_{n=-\infty}^\infty\frac{(-1)^n q^{5n^2}}{1+q^{10n-4}}
\nonumber\\[3pt]
&=\frac{2J_{10}^4J_{5,10}}{J_{1,2}^6}B_1(q)
+2q^{-1}m(-q,q^{10},q^5),
\end{align}
where $B_1(q)$ is defined as in \eqref{10-eq-B1}.
When $q$ is replaced by $-q$, \eqref{a3} becomes
\begin{align}\label{phi-q}
\phi(-q)&=q^{-1}m(-q,q^{10},-q)+q^{-1}m(-q,q^{10},q^2).
\end{align}
Replacing $q$, $z_0$, $z_1$ and $x$ by $q^{10}$, $-q$, $q^5$ and $-q$ in \eqref{mm}, respectively, we get
\begin{equation}\label{mm-3}
  m(-q,q^{10},q^{5})-m(-q,q^{10},-q)
  =-\frac{qJ_{10}^3J_{3,10}\overline{J}_{4,10}}{J_{2,10}J_{5,10}\overline{J}_{1,10}\overline{J}_{4,10}}.
\end{equation}
Replacing $q$, $z_0$, $z_1$ and $x$ by $q^{10}$, $q^2$, $q^{5}$ and $-q$ in \eqref{mm}, we obtain that
\begin{equation}\label{mm-4}
  m(-q,q^{10},q^{5})-m(-q,q^{10},q^2)
  =\frac{q^2J_{10}^3J_{3,10}\overline{J}_{2,10}}
  {J_{2,10}J_{5,10}\overline{J}_{3,10}\overline{J}_{4,10}}.
\end{equation}
It follows from \eqref{mm-3} and \eqref{mm-4} that
\begin{align}\label{mock-2---3}
m(-q,q^{10},q^{5})&=\frac{1}{2}m(-q,q^{10},-q)+\frac{1}{2}m(-q,q^{10},q^2)
\nonumber\\[3pt]
&\quad-\frac{qJ_{10}^3J_{3,10}\overline{J}_{4,10}}{2J_{2,10}J_{5,10}\overline{J}_{1,10}\overline{J}_{4,10}}
+\frac{q^2J_{10}^3J_{3,10}\overline{J}_{2,10}}
  {2J_{2,10}J_{5,10}\overline{J}_{3,10}\overline{J}_{4,10}}.
\end{align}
In view of \eqref{phi-q} and \eqref{mock-2---3}, $m(-q,q^{10},q^{5})$ can be rewritten as
\begin{align}\label{mock-2---1}
m(-q,q^{10},q^{5})&=\frac{q}{2}\phi(-q)
-\frac{qJ_{10}^3J_{3,10}\overline{J}_{4,10}}{2J_{2,10}J_{5,10}\overline{J}_{1,10}\overline{J}_{4,10}}
+\frac{q^2J_{10}^3J_{3,10}\overline{J}_{2,10}}
  {2J_{2,10}J_{5,10}\overline{J}_{3,10}\overline{J}_{4,10}}.
\end{align}
Substituting \eqref{mock-2---1} into \eqref{main-result-10-2-tempt}, we obtain \eqref{mock-2-2}.

Finally, we prove \eqref{mock-2-3}.  By
 Theorem \ref{main-result-10-2-M2}, we see that \begin{align}\label{main-result-10-3-tempt}
\overline{R}_{1,3}(4,10)
&=\frac{2J_{10}^4J_{5,10}}{J_{1,2}^6}B_4(q)
-\frac{2q^4}{J_{5,10}}\sum_{n=-\infty}^\infty\frac{(-1)^n q^{5n^2+10n}}{1+q^{10n+8}}
\nonumber\\[5pt]
&=\frac{2J_{10}^4J_{5,10}}{J_{1,2}^6}B_4(q)
+\frac{2}{qJ_{5,10}}\sum_{n=-\infty}^\infty\frac{(-1)^n q^{5n^2}}{1+q^{10n-2}}
\nonumber\\[5pt]
&=\frac{2J_{10}^4J_{5,10}}{J_{1,2}^6}B_4(q)
+2q^{-1}m(-q^3,q^{10},q^5),
\end{align}
where $B_4(q)$ is given in \eqref{10-eq-B4}.
Substituting \eqref{mock-1---2} into \eqref{main-result-10-3-tempt},
we arrive at \eqref{mock-2-3}.
This completes the proof. \qed

 \vskip 0.2cm
\noindent{\bf Acknowledgments.} This work
was supported by   the National Science Foundation of China.

\end{document}